\documentclass[10pt,a4paper,reqno,final]{amsart}
\usepackage[utf8]{inputenc}
\usepackage[lmargin=2.5cm,tmargin=3cm,bmargin=3cm,rmargin=2.5cm]{geometry}

\usepackage{amsfonts}
\usepackage{mathtools}
\usepackage{amssymb}
\usepackage{fix-cm}
\usepackage{amscd,amssymb,stmaryrd}
\usepackage{amsmath}
\usepackage{graphicx}
\usepackage{subcaption}
\usepackage{fancybox}
\usepackage{amscd,amstext}

\usepackage{color}
\usepackage{epstopdf}
\usepackage{caption}
\usepackage{multicol}
\usepackage{cite}
\usepackage{amsthm}
\usepackage{bigints}
\usepackage{amsrefs}
\usepackage{showlabels} %When the paper is "final" the label will be disappeared.
\usepackage{enumitem}
\usepackage{mathrsfs}
\usepackage{tikz}
\usepackage{empheq}
\usepackage{framed}
\usepackage[normalem]{ulem}
\usepackage{nameref}
\usepackage{esint}
\usepackage[most]{tcolorbox}
\usepackage{hyperref}

\numberwithin{equation}{section}

\newcommand{\ds}{\displaystyle}

\newtheorem{theorem}{Theorem}[section]
\newtheorem{lemma}[theorem]{Lemma}

\newtheorem{definition}[theorem]{Definition}

\newtheorem{example}{Example}
\newtheorem*{theorem*}{Theorem}
\newtheorem*{lemma*}{Lemma}
\newtheorem*{conj*}{Conjecture}
\newtheorem*{corollary*}{Corollary}
\newtheorem*{proposition*}{Proposition}
\newtheorem{remark}[theorem]{Remark}
\newcommand{\rom}[1]{\uppercase\expandafter{\romannumeral #1\relax}}

\tcbset{thmstyle/.style={
    sharp corners,           % Use sharp corners for a professional look
    boxrule=0.8pt,           % Outer border thickness
    colback=white,           % White background
    colframe=black,          % Black frame color
    borderline={.4pt}{0pt}{black}, % Inner thin line (creates the double-line effect)
    fonttitle=\sffamily\scshape,
    separator sign=,
    before skip=10pt,after skip=10pt,
}}
\newtcbtheorem[number within=section, use counter=theorem]{Theorem}{Theorem}{thmstyle}{th}
\newtcbtheorem[number within=section, use counter=theorem]{Lemma}{Lemma}{thmstyle}{lem}
\newtcbtheorem[number within=section, use counter=theorem]{Proposition}{Proposition}{thmstyle}{prop}
\newtcbtheorem[number within=section, use counter=theorem]{Corollary}{Corollary}{thmstyle}{cor}
\tcbset{defstyle/.style={
    sharp corners,
    boxrule=0.5pt,             % **Thick border** for definition
    colback=white,
    colframe=black,
    fonttitle=\sffamily\scshape,
    before skip=8pt,after skip=8pt,
}}
\newtcbtheorem[number within=section, use counter=theorem]{Definition}{Definition}{defstyle}{def}
\tcbset{remstyle/.style={
    sharp corners,
    boxrule=0.5pt,
    colback=gray!10,         % **Very light gray background**
    colframe=black,
    fonttitle=\sffamily\scshape,
    before skip=8pt,after skip=8pt,
}}
\newtcbtheorem[number within=section, use counter=theorem]{Remark}{Remark}{remstyle}{rem}

\newcommand{\ol}{\overline}

\newcommand{\lan}{\langle}
\newcommand{\ran}{\rangle}

\newcommand{\R}{\mathbb{R}}
\newcommand{\C}{\mathbb{C}}

\newcommand{\N}{\mathbb{N}}

\newcommand{\cM}{\mathcal{M}}
\newcommand{\cS}{\mathcal{S}}
\newcommand{\cF}{\mathcal{F}}

\newcommand{\cD}{\mathcal{D}}

\newcommand{\sP}{\mathscr{P}}
\newcommand{\sH}{\mathscr{H}}

\newcommand{\Ga}{\alpha}
\newcommand{\Gb}{\beta}
\newcommand{\Gd}{\delta}
\newcommand{\Ge}{\varepsilon}
\newcommand{\Gg}{\gamma}

\newcommand{\Gl}{\lambda}
\newcommand{\GL}{\Lambda}
\newcommand{\Go}{\omega}

\newcommand{\Gr}{\rho}
\newcommand{\Gs}{\sigma}

\newcommand{\Gth}{\theta}

\newcommand{\Gz}{\zeta}

\DeclareMathOperator*{\supp}{supp}

\title[Approximate cancellation on twisted Hardy spaces]{Approximate cancellation, singular integral operators, and duality on twisted Hardy spaces}
\author[C. H. Lau]{Chun Ho Lau}
\address{Department of Mathematics, National Taiwan Normal University, No. 88, Sec.4, Ting-Chou Road, Taipei, 116059, Taiwan}
\email{lauchho@ntnu.edu.tw}
\thanks{The author was supported by the National Science and Technology Council of Taiwan under research grant number 115-2811-M-003-007.} 
\date{30 November 2026}
\begin{document}

\begin{abstract}
We establish atomic and molecular characterizations of $\mathscr{H}^p(\mathbb{C}^n)$ involving approximate cancellation conditions for $0<p\leq 1$. As an application of the decomposition, we obtain boundedness results for twisted convolution operators as well as non-convolution-type (strongly) singular integral operators under suitable assumptions on their kernels and their adjoints. We also establish the duality of $\mathscr{H}^p(\mathbb{C}^n)$ using the twisted derivative and twisted Campanato space.  
\end{abstract}

\maketitle

\section{Introduction}
\subsection{Hardy spaces and local Hardy spaces}
The theory of Hardy spaces on $\R^n$ is well established and underwent substantial development during the 1970s and 1980s; see, for example, \cites{FeffStein, AM, Coifman1, Latter, TW}. A tempered distribution $f$ is said to belong to $H^p(\R^n)$ if and only if there exists a suitable class $\mathcal{F}$ of Schwartz functions defined by finitely many quasinorm bounds such that 
$$\sup_{\phi\in\cF} \sup_{t>0}|f\ast \phi_t(x)|\in L^p(\R^n),$$
where $\phi_t(x)=t^{-n}\phi(x/t).$
Real Hardy spaces provide a natural substitute for Lebesgue spaces in the range $0<p\leq 1$. Moreover, in contrast to $L^1(\R^n)$, $H^1(\R^n)$ admits a Banach space predual. Among the various characterizations of Hardy spaces, we recall those relevant to the present work.
For each $f\in H^p(\R^n)$, there exists a sequence of $L^q$ ($1<q\leq \infty$) functions $\{a_j\}_{j}$ and a sequence of scalar $\{\Gl_j\}_j\in \ell^p$ where $a_j$ is compactly supported in a ball $B(c_j,r_j)$ with $\|a\|_{L^q}\leq |B_j|^{\frac{1}{q}-\frac{1}{p}}$ and $\int a_j (x-c_j)^{\Ga}=0$ for all $|\Ga|\leq \lfloor n(1/p-1)\rfloor$, and $f=\sum_j\Gl_ja_j$ in $\cS'(\R^n)$ and $H^p(\R^n)$. Moreover, $$\|f\|_{H^p}\simeq \inf_{f=\sum_j\lambda_ja_j}\left(\sum_j|\lambda_j|^p\right)^{1/p}.
$$
This characterization is called atomic decomposition, and each $a_j$ is called an $H^p(\R^n)$ atom. Another characterization is called molecular decomposition. More precisely, there exists a sequence of functions $\{M_j\}_j$ (called molecules) and a sequence of scalars $\{\Gl_j\}_j\in \ell^p$ such that $M_j$ satisfying $\|M_j\|_{L^q(\R^n)}+\|M_j(|\cdot-c_j|/r_j)^{\Gl/q}\|_{L^q(\R^n)}\leq Cr_j^{n(\frac{1}{q}-\frac{1}{p})}$ (with suitable $\Gl$),  $\int M_j (x-c_j)^{\Ga}=0$ for all $|\Ga|\leq \lfloor n(1/p-1)\rfloor$, and $f=\sum_j\Gl_jM_j$ in $\cS'(\R^n)$ and $H^p(\R^n)$. 

These decompositions play a fundamental role in the study of boundedness of operators; if $T:\cS'(\R^n)\to\cS'(\R^n)$ is linear and continuous on $L^2(\R^n)$, then one shows that $T$ can be extended to a bounded operator on $H^p(\R^n)$ by showing $\|Ta\|_{H^p}\leq C$ for all atoms $a$ due to the equivalence between the finite and infinite atomic quasi-norms; see \cites{Bownik, MSV1, MSV2}.  For suitable singular integral operators, this can be achieved by showing that $Ta$ is an $H^p$ molecule. If $a$ is an $H^p$-atom supported in $B(c_B,r)$, then formally,
$$\int_{\R^n}Ta(x)(x-c_B)^{\Ga}dx=\lan a, T^*[(\cdot -c_B)^{\Ga}]\ran.$$ Consequently, the cancellation of $Ta$ is closely related to the action of $T^*$ on polynomials. For a sufficiently regular convolution operator, one has 
$$T^*(\mathscr{P}_{N_p})\subset \mathscr{P}_{N_p},$$
where $\mathscr{P}_{N_p}$ denotes the space of polynomials with degree at most $N_p=\lfloor n(\frac{1}{p}-1)\rfloor$. However, for non-convolution-type operators, this property does not generally hold. 

The atomic decomposition of Hardy spaces shows that the global cancellation is necessary, and therefore, Hardy spaces are not closed under multiplication by smooth cutoff functions.  Motivated by this limitation, Goldberg \cite{Gold1} introduced the local Hardy spaces, $h^p(\R^n)$ for $0<p\leq \infty$. It is another suitable generalization of $L^p(\R^n)$ for $p>1$ and $H^p(\R^n)\hookrightarrow h^p(\R^n)$ for all $0<p\leq1$  and $H^1(\R^n)\hookrightarrow L^1(\R^n)$ and indeed the inclusion is strict. We say $f\in h^p(\R^n)$ iff there exists a finite family of seminorms of Schwartz class, $\cF$, such that 
$$\sup_{\phi\in\cF} \sup_{t\in (0,1)}|f\ast \phi_t(x)|\in L^p(\R^n).$$ 
Goldberg also established an atomic characterization in which cancellation is required only for atoms supported on small balls. Later, Dafni \cite{Dafni1} showed that one can relax the cancellation condition by 
$$\bigg|\int a(x) (x-c_B)^{\Ga}dx\bigg|\leq C \min\{r^{\eta},1\} \text{ for all }|\Ga|\leq \lfloor n(1/p - 1)\rfloor, $$
where $\supp(a)\subset B(c_B,r)$. In the case $p=1$, Dafni and Yue \cite{DY} introduced a weaker (logarithmic) cancellation condition for atoms, namely,
$$\bigg|\int a(x) dx\bigg|\leq[\log(1+C/r)]^{-1}.$$
The molecular theory for $h^p$ was introduced by Komori \cite{Komori} for $\frac{n}{n+1}<p<1$. His definition is analogous to that of the $H^p$ molecule with the cancellation condition $|\int M|\leq C$. Later, Dafni and Liflyand \cite{DL} established the molecular theory for $h^1(\R)$ with the logarithmic cancellation condition. More recently, the approximate atomic and molecular characterizations were established for $h^p(\R^n)$ for all $0<p\leq 1$ in \cite{DLPV1}.

These approximate cancellation conditions allow us to impose a weaker condition on $T^*[(\cdot-c)^{\Ga}]$. It is not immediate to verify whether $T^*[(\cdot-c)^{\Ga}]=0$ on small balls. Moreover, it allows us to check whether a function is in $h^p(\R^n)$ easily. Moreover, the approximate cancellation conditions reveal the special role of the critical exponents $p=\frac{n}{n+k}$, where $k\in\N\cup\{0\}$.

\subsection{Twsited convolution and twisted Hardy spaces}
A different convolution structure arises naturally on $\C^n$; one can consider an alternative convolution operation called twisted convolution, given by
$$f\times g(z)=\int_{\C^n} f(w)g(z-w)e^{\frac{i}{2}\Im(z\cdot \ol{w})}dw.$$
The twisted convolution is closely related to the Heisenberg group and the Weyl transform. The group $\mathbb{H}^n$ are elements in $\R^n\times \R^n\times \R$ with the group action $[(x,y,t),(x',y',t')]=(x+x',y+y',t+t'+\frac{x'\cdot y-x\cdot y'}{2})$. We write the one-variable inverse Fourier transform of a function defined on $\mathbb{H}^n$ to be $F^{\Gl}(z)$ (here $z\in \C^n\cong \R^{2n}$). Then, 
$$(F\ast G)^{\Gl}(z)=\int_{\C^n} F^{\Gl}(w) G^{\Gl}(z-w)e^{\frac{i}{2}\Gl \Im(z\cdot\ol{w})}dw;$$
in particular, when $\Gl=1$, this corresponds to the twisted convolution. We also have an analogue of the Fourier transform for the convolution identity for the Weyl transform, $W(f\times g)=W(f)W(g)$. For a detailed discussion on Fourier analysis on the Heisenberg group and the Laplacian on the Heisenberg group, we refer the reader to \cites{Tbook, Stein1}.

This convolution structure naturally leads to defining the corresponding Hardy spaces, which we denote by $\sH^p$. More precisely, 
if we define $\cM_{\cF}(f)(z):=\sup_{\phi\in \cF}\sup_{t>0}|f\times \phi_t(z)|$ and $\cM_{\cF,1}(f)(z):=\sup_{\phi\in \cF}\sup_{0<t<1}|f\times \phi_t(z)|$. The tempered distribution $f$ is an element in $\sH^p$ if $\|\cM_{\cF}(f)\|_{L^p}<\infty$. It turns out 
$$\|\cM_{\cF}(f)\|_{L^p}\simeq \|\cM_{\cF,1}(f)\|_{L^p}.$$
As in the classical setting, the class $\cF$ can be chosen to be independent of $f$. Therefore, there is no difference between ``twisted Hardy spaces" and ``twisted local Hardy spaces". Indeed, Mauceri, Picardello, and Ricci \cite{MPR} firstly introduced the notion of twisted Hardy spaces. They established an atomic characterization, a characterization in terms of the Riesz transform for $\sH^1(\C^n)$, and proved $(\sH^1)^*=\mathscr{BMO}$, the twisted bounded mean oscillation space. Huang and Xing \cite{HX}, and Huang \cite{H} proved other characterizations of $\sH^1$. Huang and Wang \cite{HW} gave an atomic decomposition for $\sH^p$ for $\frac{2n}{2n+1}<p<1$ and established the boundedness of the Weyl multiplier on $\sH^p$. For the same range of $p$, Huang and Liu \cite{HL} gave a molecular decomposition on $\sH^p$ (including the endpoint $p=1$) and established the dual space of $\sH^p$. They also obtained boundedness results for local Riesz transforms as an application of their molecular characterization. Recently, Basak and Jotsaroop established an atomic decomposition for $\sH^p$ for $p\leq \frac{2n}{2n+1}$. They identified these spaces with the Hardy spaces associated with the twisted Laplacian and obtained sharp estimates for the corresponding wave operator.

%by Huang and Wang for $\frac{2n}{2n+1}<p<1$ \cite{HW}, and Basak and Jotsaroop \cite{BJ} for $0<p\leq \frac{2n}{2n+1}$. In \cite{MPR}, they have established an atomic decomposition for $\sH^p$ and showed that the dual space of $\sH^p$ is the twisted BMO (denoted by $\mathscr{BMO}$), as well as the characterization of $\sH^p$ using the twisted Riesz transform. Other characterizations of $\sH^1$ can be found in \cite{HX} and \cite{H}. In \cite{HW}, they gave an atomic decomposition for the range $\frac{2n}{2n+1}<p<1$ and proved the boundedness of the Weyl multiplier on $\sH^p$. For this range of $p$, Huang and Liu \cite{HL} gave a molecular decomposition on $\sH^p$ including $p=1$.  They used this to study the Hardy spaces associated with the twisted Laplacian. 
 \subsection{Objective of this paper}
The first goal of this paper is to establish atomic and molecular characterizations of $\sH^p$ involving approximate cancellation conditions for all $0<p\leq 1$. Motivated by the results for local Hardy spaces in \cite{DLPV1}, the exact cancellation conditions of twisted atoms and molecules can be replaced by weaker cancellation conditions.

More precisely, let $\Gg_p:=2n(\frac{1}{p}-1)$ and $N_p=\lfloor \Gg_p\rfloor$. For an atom or a molecule $H$ associated with a ball $B(c,r)$ with $r<1$, the cancellation conditions can be replaced by 
$$\bigg|\int_{\C^n}H(z)(z-c)^{\Ga}\ol{(z-c)}^{\Gb} e^{\frac{i}{2}\Im(c\cdot\ol{z})}dz\bigg|\lesssim \begin{cases}
1, & |\Ga|+|\Gb|<\Gg_p,\\
[\log(1+r^{-1})]^{-1/p},
& |\Ga|+|\Gb|=\Gg_p.
\end{cases}$$
In particular, a logarithmic cancellation condition
appears when $\Gg_p=N_p\in \N$.

As an application of these characterizations, we study twisted singular integral operators as well as non-convolution-type (strongly) singular integral operators. Different from the classical convolution setting, approximate cancellation conditions are particularly useful to establish the boundedness of a twisted singular integral operator on $\sH^p$, since the action of the adjoint of a twisted singular integral operator on $e^{\frac{i}{2}\Im(z\cdot \ol{c})}$ is not necessarily a constant multiple of the original function even if the kernel has vanishing principal-value integral. Consequently, the image of a twisted atom under such an operator does not necessarily satisfy the exact cancellation conditions of a molecule. This is the place where approximate molecular cancellation becomes useful. 

Under suitable size, smoothness, and cancellation conditions on the kernel and the adjoint of the operator, we obtain boundedness results for twisted convolution operators on $\sH^p$. We also extend these results to non-convolution-type singular integral operators and strongly singular integral operators.

The second goal of this paper is to establish the duality of $\sH^p$ for $0<p\leq \frac{2n}{2n+1}$. We first establish 
$$(\sH^{p}(\C^n))^*\cong \GL^\times_{\Gg_p}\quad \text{for }0<p<1,$$
where $\GL^\times_{\Gg_p}$ denotes a twisted Campanato space.

We then characterize these spaces using twisted
derivatives and twisted finite differences. Unlike the classical derivative, the twisted derivative does not reduce the degree of a polynomial. Thus, the classical argument based on Taylor polynomials cannot be applied directly. Nevertheless, we show that for every $\Ga>0$,
$$
    \GL^\times_{\Ga}
    \cong \widetilde{\GL}^{\times}_{\Ga}
    \cong \widetilde{\GL'}^{\times}_{\Ga},$$
with equivalent norms, where the latter two spaces are defined in terms of twisted derivatives and twisted finite differences, respectively.

The structure of this paper is as follows: we introduce the definitions of twisted convolution, atoms, molecules, and Taylor expansion in Section \ref{sec2}. In Section \ref{sec3}, we give a proof of atomic and molecular decomposition for $\sH^p$. In Section \ref{sec6}, we prove the equivalence between finite and infinite atomic decompositions. We discuss the twisted convolution singular integral operators and their non-convolution version, and then prove their boundedness on $\sH^p$ in Section \ref{sec4}. In Section \ref{sec5}, we establish the dual space of $\sH^p$, and establish the equivalence definitions using the notion of Zygmund spaces.

\section{Preliminary} \label{sec2}
\subsection{Twisted Hardy spaces}
We start with the definitions of grand maximal functions.
\begin{definition}
    Let 
    $$\cF_N=\bigg\{\phi\in C_c^\infty(B(0,1)):\sup_{|\alpha|\leq N}\|\partial^\alpha\phi\|_\infty\leq 1 \bigg\}.$$
    We define 
    $$\cM^{\times}_{\Gs}f(z):=\sup_{\phi\in \cF_N}\sup_{0<t<\Gs}|f\times \phi_t(z)|,$$
     where $\times$ denotes the twisted convolution operator given by
    $$f\times g(z):=\int_{\C^n}f(w)g(z-w)\Go(z,w)dw=\int_{\C^n}f(z-w)g(w)\Go(w,z)dw$$
    and $\Go(z,w)$ is given by $\Go(z,w):=e^{\frac{i}{2}\Im(z\cdot \overline{w})}$.
\end{definition}

To state the equivalence of twisted Hardy spaces, we need to introduce the notion of atoms. Following the idea in \cite{DLPV1}, we also introduce the notion of approximate atoms. Throughout this paper, we write $\Gg_p:=2n(\frac{1}{p}-1)$ and $N_p:=\lfloor \Gg_p\rfloor$

\begin{definition}
Let $0<p\leq 1\leq s\leq \infty$ with $p<s$.
\begin{itemize}
    \item For fixed $C>0$, we say a measurable function $A$ is a $(p,s)$-atom if there exists $B\subset \C^n$ such that 
		\begin{itemize} \item[\textnormal{(A1)}] $\supp(A)\subset B$;
			\item[\textnormal{(A2)}] $\|A\|_{L^s(\C^n)}\leq [r(B)]^{2n(1/s-1/p)}$;
			\item[\textnormal{(A3)}] $\displaystyle \int_{\C^n} A(z)(z-z_B)^{\alpha}\ol{(z-z_B)}^{\Gb}\Go(z_B,z)dz=0$ if $r(B)<1$ and $|\Ga|+|\Gb|\leq N_p$.
		\end{itemize}
    \item  We say a measurable function $A$ is an approximate $(p,s)$-atom if $A$ satisfies (A1), (A2) above, and 
		\begin{itemize}
			\item[\textnormal{(A3')}] $\displaystyle \left\{ \begin{array}{ll} \displaystyle \left| \int_{\C^n}{A(z) (z-z_B)^{\alpha}\ol{(z-z_B)}^{\Gb}\Go(z_B,z)dz} \right| \leq C, &\quad \text{if } |\alpha|+|\Gb|<\gamma_p, \\ 
				& \\
				\displaystyle \left| \int_{\C^n}{A(z) (z-z_B)^{\alpha}\ol{(z-z_B)}^{\Gb} \Go(z_B,z)dz}  \right| \leq \frac{C}{[\ln(1+\frac{1}{r(B)})]^{\frac{1}{p}}}&\quad \text{if } |\alpha|+|\Gb|=N_p =\gamma_p.
			\end{array} \right.$ 
		\end{itemize}
\end{itemize}
\end{definition}
\begin{remark}\label{rem:HigherCancellation}
As in the classical case (see \cite{Stein1}*{Chapter III, Section 2.4.2}), the atomic characterization of $\sH^p$ remains valid if we require the atoms to satisfy the cancellation conditions in (A3) up to any integer order $N\geq N_p$, and $N$ depends on $p$. Moreover, for $p<q<\infty$, the corresponding finite atomic quasi-norms are equivalent. See Section \ref{sec6}.
\end{remark}

\begin{definition}
For $0<p\leq 1\leq q\leq \infty$ with $p\neq q$, we write 
$$\sH^{p,q}_{fin}:=\{f\in \cS'(\C^n): f=\sum_{j=1}^{N}\Gl_ja_j, \ a_j \text{ are $(p,q)$-atoms}.\}$$
with 
$$\|f\|_{\sH^{p,q}_{fin}}:=\inf\bigg\{\bigg(\sum_{j=1}^N|\Gl_j|^p\bigg)^{\frac{1}{p}}\bigg\},$$
where the infimum is taken over all such finite atomic decompositions.
We also write 
$$\sH^{p,q}_{at}:=\{f\in \cS'(\C^n): f=\sum_{j=1}^{\infty}\Gl_ja_j, \ a_j \text{ are $(p,q)$-atoms}.\}$$
and 
$$\|f\|_{\sH^{p}_{at}}:=\inf\bigg\{\bigg(\sum_{j=1}^\infty|\Gl_j|^p\bigg)^{\frac{1}{p}}\bigg\},$$
where the infimum is taken over all such atomic decompositions.
\end{definition}

We have the following.
\begin{theorem}[\cites{BJ, MPR}]
Let $0<p\leq1$ and $f\in \cS'(\C^n)$. Then, for sufficiently large $N$ depending on $p$, we have 
$$\|\cM^{\times}_1f\|_{L^p(\C^n)}\simeq \|\cM^{\times}_{\infty}f\|_{L^p(\C^n)}\simeq \|f\|_{\sH^p_{at}}.$$
In this case, we say $f\in \sH^p$ if one of the quasinorms is finite.
\end{theorem}

We also introduce the notion of molecules. 
\begin{definition}
Let $0<p\leq 1\leq s<\infty$ with $p<s$ and $\Gl>2n(\frac{s}{p}-1)$. 
\begin{itemize}
    \item  For fixed $C>0$, we say a measurable function $M$ is a $(p,s,\Gl)$-molecule if there exists $B\subset \C^n$ such that 
		\begin{itemize} \item[\textnormal{(M1)}] $\|M\|_{L^s(2B)} \leq C \, r(B)^{2n\left(\frac{1}{s}-\frac{1}{p}\right)}$
			\item[\textnormal{(M2)}] $\displaystyle \| M \, |\cdot - z_B|^{\frac{\lambda}{s}} \,  \|_{L^{s}((2B)^c)} \leq C \, r(B)^{\frac{\lambda}s + 2n \left(\frac{1}{s}-\frac{1}{p}\right)}$
			\item[\textnormal{(M3)}] $\displaystyle \int_{\C^n} M(z)(z-z_B)^{\alpha}\ol{(z-z_B)}^{\Gb}\Go(z_B,z)dz=0$ if $r(B)<1$ and $|\Ga|+|\Gb|\leq N_p$.
		\end{itemize}
    \item  We say a measurable function $M$ is an approximate $(p,s,\Gl)$-molecule if $M$ satisfies (M1), (M2) above, and 
		\begin{itemize}
			\item[\textnormal{(M3')}] $\displaystyle \left\{ \begin{array}{ll} \displaystyle \left| \int_{\C^n}{M(z) (z-z_B)^{\alpha}\ol{(z-z_B)}^{\Gb}\Go(z_B,z)dz} \right| \leq C', &\quad \text{if } |\alpha|+|\Gb|<\gamma_p, \\ 
				& \\
				\displaystyle \left| \int_{\C^n}{M(z) (z-z_B)^{\alpha}\ol{(z-z_B)}^{\Gb} \Go(z_B,z)dz}  \right|\leq \frac{C}{[\ln(1+\frac{1}{r(B)})]^{\frac{1}{p}}}&\quad \text{if } |\alpha|+|\Gb|=N_p =\gamma_p.
			\end{array} \right.$ 
		\end{itemize}
\end{itemize}
\end{definition}
\begin{remark}
     Indeed, we can also use $\C^n$ instead of $2B$ in (M1) and $(2B)^c$ in (M2), see \cite{DLPV1}*{Remark 3.6}. This has nothing to do with the structure of $\C^n$. 
     %Indeed, using (M1) and (M2), we have 
    % \begin{align*}
    %     \|M\|_{L^s(\C^n)}\leq \|M\|_{L^s(2B)}+r^{-\frac{\Gl}{s}}\|M|\cdot-z_B|^{\frac{\Gl}{s}}\|_{L^s([2B]^c)}\leq 2C r^{2n(\frac{1}{s}-\frac{1}{p})},
    % \end{align*}
    % and 
    % \begin{align*}
    %     \|M|\cdot-z_B|^{\frac{\Gl}{s}}\|_{L^s(\C^n)}\leq \|M\|_{L^s(2B)}+r^{-\frac{\Gl}{s}}\|M|\cdot-z_B|^{\frac{\Gl}{s}}\|_{L^s([2B]^c)}\leq 2C r^{\frac{\Gl}{s}+2n(\frac{1}{s}-\frac{1}{p})}.
    % \end{align*}
\end{remark}
We note that every atom is a $(p,s,\Gl)$-molecule for any $\Gl>2n(s/p-1)$.
\begin{remark}
We remark that (A3') or (M3') holds for free whenever $r(B)\geq 1$. Indeed, when $r(B)\geq1$, the conditions (A3') and (M3') follow from the corresponding size conditions. By (M1) and (M2) (or (A1) and (A2) for atoms), we have
    $$\int_{\C^n}|M(z)||z-z_B|^{|\Ga|+|\Gb|}dz\lesssim r(B)^{|\Ga|+|\Gb|-\Gg_p}\leq 1$$
    for all $|\Ga|+|\Gb|\leq N_p$. Moreover, $[\ln(1+r(B)^{-1})]^{-1/p}\geq1$ when $r(B)\geq1$.
    Thus, no additional cancellation conditions are imposed on approximate atoms and molecules associated with balls of radius at least one.
\end{remark}

\subsection{Twisted derivatives and twisted Taylor polynomials}
We write $z=$ $(\Re(z_1),\cdots, \Re(z_n),$ $\Im(z_1),$ $\cdots, \Im(z_n))\in \R^{2n}$ and 
$$X_j:=\begin{cases}
    \frac{\partial}{\partial x_j}+i\frac{y_j}{2}, &\quad j\in\{1,\cdots,n\},\\
    \frac{\partial}{\partial y_{j-n}}-i\frac{x_{j-n}}{2}, &\quad j-n\in\{1,\cdots,n\};\\
\end{cases} \quad \widetilde{X}_j:=\begin{cases}
    \frac{\partial}{\partial x_j}-i\frac{y_j}{2}, &\quad j\in\{1,\cdots,n\},\\
    \frac{\partial}{\partial y_{j-n}}+i\frac{x_{j-n}}{2}, &\quad j-n\in\{1,\cdots,n\};\\
\end{cases}. $$
In this note, we say $X_j$ to be the left-invariant derivative and $\widetilde{X}_j$ to be the right-invariant derivative. One has the following properties:
\begin{itemize}
    \item $X_j(f\times g)= f\times (X_jg)$;
    \item $\widetilde{X}_j(f\times g)=(\widetilde{X}_jf)\times g $
\end{itemize}
\begin{remark}\label{rem:regularity}
    By the definitions of $X_j$ and $\widetilde{X}_j$, every ordinary derivative of order $k$ can be expressed as a linear combination of twisted derivatives of order at most $k$ with polynomial coefficients, and vice versa. Consequently, $g\in C^k(\C^n)$ if and only if all its distributional twisted derivatives up to order $k$, including $g$ itself, are continuous. This holds for both left- and right-invariant
    derivatives.
\end{remark}
\begin{definition}
    The left and right twisted translation of $f$ by $w$ is defined to be 
    $$\tau_w^{X}f(z)=f(z+w)\Go(z,w)$$
    and 
    $$\tau_w^{\widetilde{X}}f(z)=f(z-w)\Go(z,w),$$
    respectively.
\end{definition}
It can be seen that $[\widetilde{X}_j,\tau^{X}_w]=0$ and $[X_j,\tau^{\widetilde{X}}_w]=0$ for $j\in \{1,\cdots,2n\}.$ Also, the translations are not commutative, because 
$$\tau^{\widetilde{X}}_{a}\tau^{\widetilde{X}}_{b}f(z)=\tau^{\widetilde{X}}_a[f(\cdot-b)\Go(\cdot,b)](z)=f(z-a-b)\Go(z-a,b)\Go(z,a)=\ol{\Go(a,b)}\tau^{\widetilde{X}}_{b+a}f(z)$$
and 
$$\tau^{\widetilde{X}}_{a}\tau^{\widetilde{X}}_{b}f(z)=[\ol{\Go(a,b)}]^2\tau^{\widetilde{X}}_{b}\tau^{\widetilde{X}}_{a}f(z)$$
for arbitrary $a,b\in\C^n$. However, if $a=th$ and $b=sh$ with $s,t\in \R$ and $h\in \C^n$, then we have 
$$\tau^{\widetilde{X}}_{sh}\tau^{\widetilde{X}}_{th}=\tau^{\widetilde{X}}_{(s+t)h}=\tau^{\widetilde{X}}_{th}\tau^{\widetilde{X}}_{sh}.$$

In \cite{BJ}, Basak and Jotsaroop obtained the following Taylor expansions of degree $N$ centered at $z-z_c$ with integral error:
\begin{align}\label{rightT} %right
    &h(z-w)\Go(z-z_c,w-z_c)\nonumber\\
    &=h(z-z_c)+\sum_{k=1}^{N}\sum_{1\leq i_1,\cdots, i_k\leq 2n} \frac{(-1)^k}{k!} T_{i_1,\cdots, i_k}(w-z_c)[h(z-z_c)]\\
    &\quad +\sum_{i_1,\cdots,i_{N+1}=1}^{2n}\int_0^1\frac{(-1)^{N+1}}{N!}(1-\Gth)^{N}e^{\frac{i}{2}\Gth\Im((z-z_c)\cdot\ol{(w-z_c)})}T_{i_1,\cdots,i_{N+1}}(w-z_c)[h((z-z_c)-\Gth (w-z_c))]d\Gth\nonumber
\end{align}
and 
\begin{align}\label{leftT} %right
    &h(z-w)\Go(z-z_c,w-z_c)\nonumber\\
    &=h(-w+z_c)+\sum_{k=1}^{N}\sum_{1\leq i_1,\cdots, i_k\leq 2n} \frac{1}{k!} U_{i_1,\cdots, i_k}(z-z_c)[h(-(w-z_c))]\\
    &\quad +\sum_{i_1,\cdots,i_{N+1}=1}^{2n}\int_0^1\frac{1}{N!}(1-\Gth)^{N}e^{\frac{i}{2}\Gth\Im((z-z_c)\cdot\ol{(w-z_c)})}U_{i_1,\cdots,i_{N+1}}(z-z_c)[h(\Gth (z-z_c)-(w-z_c))]d\Gth, \nonumber
\end{align}
where $u_j$ is defined by $u_{j}=\Re(w_j)$ and $u_{j+n}=\Im(w_j)$ for $j\in\{1,\cdots, n\}$, and 
$$T_{i_1,\cdots,i_k}(w)[F(z')]=u_{i_1}\cdots u_{i_k} \big(\widetilde{X}_{i_1}\cdots\widetilde{X}_{i_k}[F]\big)\big|_{(x,y)=(\Re(z'),\Im(z'))},$$
$$U_{i_1,\cdots,i_k}(w)[F(z')]=u_{i_1}\cdots u_{i_k} \big({X}_{i_1}\cdots{X}_{i_k}[F]\big)\big|_{(x,y)=(\Re(z'),\Im(z'))}.$$

In other words, the function $h(z-w)$ can be expanded using twisted Taylor polynomials, which are a product of a polynomial $P$ centered at $z-z_c$ with the phase $\Go(\cdot-z_c, z-z_c)$.

\section{Approximate atomic and molecular decompositions}\label{sec3}
\begin{theorem}\label{atomicdecompo}
    For any $f\in \sH^p(\C^n)$, there exist families of approximate $(p,s)$-atoms $\{A_j\}_j$ and complex numbers $\{\Gl_j\}_j$ such that $\sum_j|\Gl_j|^p<\infty$ and 
    $$f=\sum_{j} \Gl_j A_j$$
    in the sense of distributions.
\end{theorem}

Note that if $f\in\sH^p(\C^n)$, then we already knew that $f$ attained an atomic decomposition using $(p,s)$-atoms. Thus, it remains to show that $f$ has a representation using approximate $(p,s)$-atoms then $f\in \sH^p(\C^n)$.

\begin{lemma}\label{atomic1}
   If $A$ is an approximate $(p,s)$-atom, then $\|A\|_{\sH^p}\lesssim 1$. 
\end{lemma}
\begin{proof}
    It suffices to show that 
    $$\int_{\C^n} \sup_{\phi\in\cF_N}\sup_{0<t<1}|A\times \phi_t(z)|^p dz\leq C,$$
  where $C$ is independent of $A$ and $N$ is sufficiently large  (at least $N_p+1$). Then, by the equivalent norms on $\sH^p$, we can conclude that $A\in \sH^p(\C^n)$. Let $A$ be an approximate $(p,s)$-atom with $\supp(A)\subset B(z_0,r)$.

    Note that 
    \begin{align*}
        |A\times \phi_t (z)| &\leq \int_{\C^n}|A(w)||\phi_t(z-w)|dw \\
        &\leq  \frac{\|\phi\|_{L^\infty}}{|B(z,t)|}\int_{B(z,t)}|A(w)|dw\leq  M_{HL}(A)(z),    \end{align*}
        where $M_{HL}$ is the Hardy-Littlewood maximal function. 
        Then, 
    $$\int_{B(z_0,2r)}\sup_{\phi\in \cF_N}\sup_{0<t<1}|A\times \phi_t (z)|^pdz \leq C_{\phi,s} |B(z_0,2r)|^{1-\frac{p}{s}}\|A\|_{L^s}^{p}\leq C_{s,p,n}$$
    using (A2) provided that $p\in (0,1]$ and $s>1$. If $s=1$ (and $p<1$), then one can use the weak $(1,1)$ estimate for $M_{HL}$ to obtain the desired estimate, see \cite{WeightedIneq}*{Lemma 3.1}.

    On the other hand, by \eqref{rightT}
    \begin{align*}
       \phi_t(z-w)\omega(z,w)&= \phi_t((z-z_0)-(w-z_0))\omega(z-z_0, w-z_0) \omega(z,z_0)\omega(z_0,w) \\
    &=:I(z,w)+II(z,w),
    \end{align*}
    where 
    \begin{align*}
       I(z,w):=\phi_t(z-z_0)\omega(z,z_0)\omega(z_0,w) +\sum_{k=1}^{N_p}\sum_{1\leq i_1,\cdots, i_k\leq 2n} \frac{(-1)^k}{k!} T_{i_1,\cdots, i_k}(w-z_0)[ \phi_t(z-z_0)]\omega(z,z_0)\omega(z_0,w) 
    \end{align*}
    and
    \begin{align*}
        II(z,w)&=\omega(z,z_0)\omega(z_0,w) \bigg( \sum_{i_1,\cdots,i_{N_p+1}=1}^{2n}\int_0^1\frac{(-1)^{N_p+1}}{N_p!}(1-\Gth)^{N_p}e^{\frac{i}{2}\Gth\Im((z-z_0)\cdot\ol{(w-z_0)})}\\
        &\quad\quad\quad\quad\quad\quad\quad\quad\quad\quad\quad\quad\quad\quad \quad T_{i_1,\cdots,i_{N_p+1}}(w-z_0)[  \phi_t((z-z_0)-\Gth (w-z_0))]d\Gth\bigg) \\
    \end{align*}
    Then,
    $$|A\times \phi_t(z)|\leq \bigg|\int_{\C^n} A(w)I(z,w)dw\bigg|+\int_{\C^n}|A(w)||II(z,w)|dw.$$

    We first consider the term with $I(z,w)$. Note that $\omega(z_0,z)$ is independent of $w$, we can drop the factor $\omega(z_0,z)$. Moreover, since, after expanding and dropping the factor, one can rewrite $I$ as 
    $$I(z,w)=\omega(z_0,w)\sum_{|\Ga|+|\Gb|\leq N_p}K_{\Ga,\Gb,t}(z-z_0)(w-z_0)^{\Ga}\ol{(w-z_0)}^{\Gb},$$
    where $K_{\Ga,t}(z-z_0)$ depends only on $\phi_t$ and its partial derivatives up to order $|\Ga|$. Note that there are two types of terms: one only involves derivatives of $\phi_t$, and one is of the form $P(z-z_0)\partial^{\Gb}\phi_t(z-z_0)$. If $|z-z_0|\in [2r,2]$, we can control as follows: $|\partial^{\Ga}\phi_t(z)|\leq |z-z_0|^{-2n-|\Ga|}$, and $|P(z-z_0)\partial^{\Gb}\phi_t(z-z_0)|\leq |z-z_0|^{\deg(P)-2n-|\Gb|}$. However, since $|z-z_0|\leq 2$ and $|\Ga|+|\Gb|\leq N_p$, we have
    $$|K_{\Ga,\Gb,t}(z-z_0)|\lesssim |z-z_0|^{-2n-|\Ga|-|\Gb|}.$$ If $|z-z_0|\geq 2$, then one can use the fact that $|\partial^{\Gb}\phi_t(z-z_0)|\leq Ct^{|\Ga|+|\Gb|+\deg(P)+\Ge}|z-z_0|^{-2n-|\Ga|-|\Gb|-\deg(P)-\Ge}$ and $t\in (0,1)$, we have $$|K_{\Ga,\Gb,t}(z-z_0)|\lesssim |z-z_0|^{-2n-|\Ga|-|\Gb|-\Ge}.$$ Here all the implicit constants are independent of $z-z_0$, $t$, and $\phi$.
    Therefore, if $r<1$,
    \begin{align*}
        &\int_{|z-z_0|\geq 2r}\sup_{t\in (0,1)}\bigg|\int_{\C^n}A(w)I(z,w)dw\bigg|^pdz\\ &=\bigg(\int_{2r\leq |z-z_0|\leq 2}+\int_{|z-z_0|\geq 2}\bigg)\sup_{t\in (0,1)}\bigg|\int_{\C^n}A(w)I(z,w)dw\bigg|^pdz \\
        &\lesssim \sum_{|\Ga|+|\Gb|\leq N_p}\int_{2r\leq |z-z_0|\leq 2} \frac{1}{|z-z_0|^{2np+p|\Ga|+p|\Gb|}}dz \bigg|\int_{\C^n}A(w)(w-z_0)^{\Ga}\ol{(w-z_0)}^{\Gb}\omega(z_0,w)dw\bigg|^p\\
        &\quad\quad + \sum_{|\Ga|+|\Gb|\leq N_p}\int_{|z-z_0|\geq 2} \frac{1}{|z-z_0|^{2np+p|\Ga|+p|\Gb|+p\Ge}} r^{p|\Ga|+p|\Gb|}
    \end{align*}
    The second term above is always finite whatever $p$, $\Ga$, and $\Gb$ by choosing $\Ge=\Gg_p+1$; moreover, it is finite and independent of $r$ because $|z-z_0|\geq 2\geq 2r$ and the integral cancels $r^{p|\Ga|+p|\Gb|}$. The first term is slightly different. If $N_p=\Gg_p=|\Ga|+|\Gb|$, then the integral gives us $C\ln(r^{-1})$, which implies the first term is finite for any $r>0$ with the upper bound from (A3'). If $\Gg_p>N_p$ or $|\Ga|+|\Gb|\leq N_p-1$, then $2np+p|\Ga|+p|\Gb|=2n-p(\Gg_p-|\Ga|-|\Gb|)<2n$, which implies the first is finite and independent of $r$ with the aid of the first line of (A3'). 

    If $r\geq 1$, then $|z-z_0|\geq 2r\geq 2$; and therefore,
    \begin{align*}
         \sum_{|\Ga|+|\Gb|\leq N_p}\int_{|z-z_0|\geq 2r} \frac{1}{|z-z_0|^{2np+p|\Ga|+p|\Gb|+p\Ge}} r^{p|\Ga|+p|\Gb|} \leq C_{n,p}
    \end{align*}
    by the same argument.

    Therefore, $$\int_{[B(z_0,2r)]^c}\sup_{\phi\in \cF_N}\sup_{t\in (0,1)}\bigg|\int_{\C^n}A(w)I(z,w)dw\bigg|^pdz \leq C_{n,p}<\infty$$
    with $C_{n,p}$ independent of $z_0$ and $r$.
    
    To handle the term $II(z,w)$, we can handle it as in \cite{BJ}. We first consider $r<1$. For $|w-z_0|<r$,
    $$|II(z,w)|\lesssim \frac{r^{N_p+1}}{t^{2n+N_p+1}}.$$
    Therefore, if $|z-z_0|\geq 2r$, then $|z-z_0|\leq  t+r\leq t+\frac{|z-z_0|}{2}$ and this implies $|z-z_0|\leq 2t$. Thus,
    $$\int_{\C^n}|A(w)||II(z,w)|dw \lesssim \frac{r^{N_p+1+2n-\frac{2n}{p}}}{t^{2n+N_p+1}}\lesssim \frac{r^{N_p+1+2n-\frac{2n}{p}}}{|z-z_0|^{2n+N_p+1}}.$$

    Thus, 
    \begin{align*}
\int_{[B(x_0,2r)]^c}\sup_{\phi\in \cF_N}\sup_{0<t<1}\bigg|\int_{\C^n}A(w)II(z,w)dw\bigg|^pdz&\lesssim \int_{[B(x_0,2r)]^c}  \frac{r^{pN_p+p+2np-2n}}{|z-z_0|^{2np+N_pp+p}}dz\leq C_{n,p}
    \end{align*}
    provided that $2np+N_pp+p>2n$, which is true as $N_p+1> 2n(\frac{1}{p}-1)$. Moreover, $C_{n,p}$ is independent of $r$.

 If $r\geq1$, then from $A\times\phi_t$, $|z-z_0-\Gth(w-z_0)|\leq t<1$ and $|w-z_0|\leq r$, which implies $A\times\phi_t$ is non-zero whenever
    $$|z-z_0|<t+r<2r. $$
    Thus, on $z\in [B(z_0,2r)]^c$, the expression $A\times \phi_t$ is indeed 0.
\end{proof}

We now establish the decomposition using approximate molecules.

\begin{theorem}\label{moleculardecompo}
    For any $f\in \sH^p(\C^n)$, there exist families of approximate $(p,s,\Gl)$-molecules $\{M_j\}_j$ and complex numbers $\{\Gl_j\}_j$ such that $\sum_j|\Gl_j|^p<\infty$ and 
    $$f=\sum_{j} \Gl_j M_j$$
    in the sense of distributions. Moreover, 
    $$\inf \bigg\{\sum_{j=1}^{\infty}|\Gl_j|^p\bigg\}\simeq \|f\|_{\sH^p}^p,$$
    where the infimum is taken over all possible approximate molecular decompositions of $f$. This is also true for $(p,s,\Gl)$-molecules.
\end{theorem}

Since $(p,s)$-atoms are $(p,s,\Gl)$-molecules and $(p,s,\Gl)$-molecules are also approximate $(p,s,\Gl)$-molecules, the atomic decomposition of $f$ automatically implies the molecular decomposition. Thus, it reduces to estimating $\|M\|_{\sH^p}\lesssim 1$. 

\begin{lemma}\label{mole1}
   If $M$ is an approximate $(p,s,\Gl)$-molecules, then $\|M\|_{\sH^p}\lesssim 1$. 
\end{lemma}
\begin{proof}
The proof is based on the proof given in \cite{WeightedIneq}*{Theorem III.7.16} and \cite{DLPV1}. Let $M$ be a normalized molecule associated to $B=B(z_B,r)\subset \C^n$. Define $B_0=E_0=B$, $B_k=B(z_B,2^kr)$, and $E_k=B_{k}\setminus B_{k-1}$. Choose the usual polynomials $\psi^k_{\Ga,\Gb}$ such that 
$$\fint_{E_k}\psi^k_{\Ga,\Gb}(z)(z-z_B)^{\Ga'}\ol{(z-z_B)}^{\Gb'}dz=1$$
if and only if $(\Ga',\Gb')=(\Ga,\Gb)$, and 0 otherwise. Set $\phi^k_{\Ga,\Gb}(z)=\psi^k_{\Ga,\Gb}(z)\Go(z,z_B)\chi_{E_k}(z)$. We further define $M_k(z)=M(z)\chi_{E_k}(z)$ and $P_k=\ds \sum_{|\Ga|+|\Gb|\leq N_p}m_{\Ga,\Gb}^k\phi_{\Ga,\Gb}^k$, where $P_k$ is a complex twisted polynomial restricted on $E_k$ of degree at most $N_p$, and satisfies 
	\begin{align}\label{AF-1'}
		 \fint_{E_k} P_k(z) (z-z_B)^{\alpha}\ol{(z-z_B)^{\Gb}}\Go(z_B,z)dz  &= m^{k}_{\alpha,\Gb}\fint_{E_k}{\phi_{\Ga,\Gb}^k(z) \  (z-z_B)^{\alpha}\ol{(z-z_B)^{\Gb}}\Go(z_B,z)dz}\nonumber\\
          &  =	m^{k}_{\alpha,\Gb}:= \fint_{E_k}{M(z) (z-z_B)^{\alpha}\ol{(z-z_B)^{\Gb}}\Go(z_B,z)dz}.
		 \end{align}
         and 
        		 \begin{equation} \label{AF-2'}
		  |P_k(z)| \leq C_{n,p}\fint_{E_k}|M(z)|dz
		   \end{equation}
		for some constant independent of $k$.   This is done by choosing the polynomials $\phi_{\alpha,\Gb}^{k}$ to have $(\gamma_1,\gamma_2)$-th moment equal to $|E_k|$ when $(\gamma_1,\gamma_2)=(\Ga,\Gb)$, and zero otherwise, and noting that $(2^kr)^{|\alpha|+|\Gb|} |\phi_{\alpha,\Gb}^k(z)|\leq C$, where $C$ is independent of $k$.

For $|\Ga|+|\Gb|\leq N_p$ and $k\in \N_0$, we set 
\begin{align*}
    \nu_{\Ga,\Gb}&=\sum_{j=0}^{\infty} |E_j|m_{\Ga,\Gb}^j=\int_{\C^n} M(z)(z-z_B)^{\Ga}\ol{(z-z_B)^{\Gb}}\Go(z_B,z)dz;\\
N_{\Ga,\Gb}^k&=\sum_{j=k+1}^{\infty}|E_j|m_{\Ga,\Gb}^j=\int_{B_k^c} M(z)(z-z_B)^{\Ga}\ol{(z-z_B)^{\Gb}}\Go(z_B,z)dz.\\
\end{align*}
Then, we have 
\begin{align*}
    \sum_{k=0}^{\infty} P_k(z)&=\sum_{|\Ga|+|\Gb|\leq N_p} \left(\sum_{k=1}^{\infty} (N_{\alpha,\Gb}^{k-1} - N_{\alpha,\Gb}^{k}) |E_k|^{-1} \phi_{\alpha,\Gb}^{k}(z) + (\nu_{\alpha,\Gb} - N_{\alpha,\Gb}^{0})|E_0|^{-1} \phi_{\alpha,\Gb}^{0}(z)\right)\\
    &=\sum_{|\Ga|+|\Gb|\leq N_p} \left(\sum_{k=0}^{\infty} N_{\alpha,\Gb}^{k} |E_{k+1}|^{-1} \phi_{\alpha,\Gb}^{k+1}(z)-\sum_{k=0}^{\infty} N_{\alpha,\Gb}^{k} |E_{k}|^{-1} \phi_{\alpha,\Gb}^{k}(z) + \nu_{\alpha,\Gb}|E_0|^{-1} \phi_{\alpha,\Gb}^{0}(z)\right)\\
	&= \sum_{k=0}^{\infty}{\bigg(\sum_{|\alpha|+|\Gb|\leq N_p}{\Phi_{\alpha,\Gb}^{k}(z)}} + \sum_{|\alpha|+|\Gb| \leq N_p}\nu_{\alpha,\Gb}\, |E_0|^{-1} \, \phi_{\alpha,\Gb}^{0}(z)\bigg),
		\end{align*}
where 
$$\Phi_{\Ga,\Gb}^{k}(z)=N_{\Ga,\Gb}^{k}[|E_{k+1}|^{-1}\phi_{\Ga,\Gb}^{k+1}(z)-|E_k|^{-1}\phi_{\Ga,\Gb}^k(z)]$$
for $k\in \N_0$. Therefore, we can decompose $M$ as 
$$M=\sum_{k=0}^{\infty} (M_k-P_k)+\sum_{k=0}^{\infty}\sum_{|\Ga|+|\Gb|\leq N_p} \Phi^k_{\Ga,\Gb}+\sum_{|\Ga|+|\Gb|\leq N_p}\nu_{\Ga,\Gb}|E_0|^{-1}\phi_{\Ga,\Gb}^0=:S_1+S_2+S_3.$$

To handle $S_1$, we first show that $M_k-P_k$ is a multiple of an $\sH^p$ atom. To see this, note that 
\begin{align*}
    \|M_k\|_{L^s}\leq 2^{\frac{\Gl}{s}}C_{n,s}|B_k|^{\frac{1}{s}-\frac{1}{p}}(2^k)^{-\frac{\Gl}{s}+2n(\frac{1}{p}-\frac{1}{s})}.
\end{align*}
and from \eqref{AF-2'} one has 
$$\|P_k\|_{L^s}\leq C_{n,p} \|M_k\|_{L^s}.$$
Moreover, from \eqref{AF-1'}, for any $|\Ga|+|\Gb|\leq N_p$ we have
$$\int_{\C^n}[M_k(z)-P_k(z)](z-z_B)^{\Ga}\ol{(z-z_B)^{\Gb}}\Go(z_B,z)dz=\int_{E_k}[M_k(z)-P_k(z)](z-z_B)^{\Ga}\ol{(z-z_B)^{\Gb}}\Go(z_B,z)dz=0.$$
This shows that $M_k-P_k$ is a constant multiple of a $(p,s)$-atom. To see it is an element in $\sH^p$, we observe that if $(M_k-P_k)(z) = t_{k} \, a_k(z)$ then $t_{k} = \| M_k-P_k \|_{L^s} \, |B_k|^{\frac{1}{p}-\frac{1}{s}}$, $a_k(z)=\dfrac{M_k(z)-P_k(z)}{\| M_k-P_k \|_{L^s}} \, |B_k|^{\frac{1}{s}-\frac{1}{p}}$, and 
		\begin{align} \label{control-t_j'}
			\sum_{k=0}^{\infty}{|t_{k}|^p} \leq 2^{\frac{\lambda p}{s}} \, C_{n,p,s} \, \sum_{k=0}^{\infty}{(2^k)^{-\frac{\lambda p}{s}+2n \left( 1-\frac{p}{s} \right)}} = C_{n,p, s,\lambda} < \infty
		\end{align}
		provided $\lambda> 2n \left( s/p-1 \right)$. We point out that the closer $\lambda$ gets to $2n(s/p-1)$, the bigger the constant appearing in \eqref{control-t_j'}. 

For $S_2$, we claim that $\Phi_{\Ga,\Gb}^j(z)$ is a constant multiple of a $(p,\infty)$ atom.  The cancellation follows from the moment conditions on $\phi_{\alpha,\beta}^{k}$.  For the size condition, from H\"older inequality and the $L^s$ norm of $M_k$, for every $|\alpha|+|\beta|\leq N_p$ and $k\in \N$, %$\lambda> n \left( {s}/{p}-1 \right)$, 
		one has
		\begin{align*} 
			|N_{\alpha,\Gb}^{k}|  \leq C_{n,p, s,\lambda} \, |B_{k}|^{1-\frac{1}{p}} (2^{k}r)^{|\alpha|+|\Gb|} (2^{k})^{-\frac{\lambda}{s}+2n \left( \frac{1}{p}-\frac{1}{s} \right)}.
		\end{align*}
		Hence, since $(2^kr)^{|\alpha|+|\Gb|} |\phi_{\alpha,\Gb}^k(x)|\leq C$ uniformly, it follows
		$$
		|N^{k}_{\alpha,\beta}|E_k|^{-1}\phi_{\alpha,\beta}^{k}(x)| \leq C_{n,p, s,\lambda} \,  |B_k|^{-\frac{1}{p}}(2^k)^{-\frac{\lambda}{s}+2n \left( \frac{1}{p}-\frac{1}{s} \right)}.
		$$
		Therefore, writing $\displaystyle{\sum_{|\alpha| +|\beta|\leq N_p} \Phi_{\alpha}^{j}(x)} =s_{k} \,\, b_{k}(x)$, where $s_{k}= C_{n,p, s,\lambda}  \, (2^k)^{-\frac{\lambda}{s}+2n \left( \frac{1}{p}-\frac{1}{s} \right)}$ 
		%$b_{k}(x)=(2^{k})^{\frac{\lambda}{n}-n \left( \frac{1}{p}-\frac{1}{s} \right)}\displaystyle{\sum_{|\alpha| \leq N_p} \Phi_{\alpha}^{j}(x)}$ 
		for some appropriate constants $C_{n,p, s,\lambda}$, we get that $b_k$ are $(p,\infty)$ atoms and		
		\begin{equation*}
			\sum_{k=0}^{\infty}{|s_{k}|^p} = C_{n,p, s,\lambda}\, \sum_{k=0}^{\infty}{(2^k)^{-\frac{\lambda p}{s}+2n \left( 1-\frac{p}{s} \right)}} < \infty,
		\end{equation*}
		where again we used that $\lambda> 2n \left( s/p-1 \right)$. 

        For $S_3$, note that 
        $$
		a_\omega = \sum_{|\alpha|+|\beta| \leq N_p}\nu_{\alpha,\beta} \, |E_0|^{-1} \, \phi_{\alpha,\beta}^{0}(x).
		$$
		This function is supported on $E_0 = B_0$, and (M1) and (M2) give
		\begin{align*}
			|\nu_{\alpha,\beta}| & = \left|\int_{\C^n}{M(z)(z-z_B)^{\alpha}dx}\ol{(z-z_B)^{\beta}}\Go(z_B,z)dz\right| \lesssim r^{|\alpha|+|\beta|+2n\left( 1-\frac{1}{p} \right)}.
			% 			& \leq \int_{|x-x_B|<2r}{|M(x)| \ |x-x_B|^{|\gamma|}dx}+\int_{|x-x_B| \geq 2r}{|M(x)| \ |x-x_B|^{|\gamma|}dx} \\
			% 			& \leq (2r)^{|\gamma|} \int_{|x-x_B|<2r}{|M(x)|dx} +\int_{|x-x_B| \geq 2r}{|M(x)| \ |x-x_B|^{\lambda -(\lambda-|\gamma|)}dx} \\
			% 			& \lesssim r^{\ |\gamma|+n\left(1-\frac{1}{p}\right)} + r^{-(\lambda-|\gamma|)}\int_{|x-x_B| \geq 2r}{|M(x)| \ |x-x_B|^{\lambda}dx} \\
			% 			& \lesssim r^{\ |\gamma|+n\left(1-\frac{1}{p}\right)} + r^{\ -(\lambda -|\gamma|)+\lambda+n\left(1-\frac{1}{p}\right) } \\
			% 			& \lesssim r^{|\gamma|+n \left(1-\frac{1}{p}\right)}.  
		\end{align*}
        The $L^s$ estimate now follows immediately from the fact that $r^{|\alpha|} |\phi_{\alpha}^0(x)|\leq C$:
		\begin{align}
			\left\| \sum_{|\alpha|+|\beta| \leq N_p}\nu_{\alpha} \ |E_0|^{-1}\phi_{\alpha,\beta}^{0} \right\|_{L^s} & \leq \sum_{|\alpha|+|\beta|\leq N_p}{|\nu_{\alpha,\beta}| \ |E_0|^{-1} \left( \int_{E_0}{|\phi_{\alpha,\beta}^{0}(x)|^sdx} \right)^{\frac{1}{s}}} \nonumber \\
			& \leq \sum_{|\alpha|+|\beta|\leq N_p}{|\nu_{\alpha,\beta}| \, |E_0|^{\frac{1}{s}-1} \ r^{-|\alpha|-|\beta|}} \lesssim \, r^{\,2n\left(\frac{1}{s}-\frac{1}{p}\right)}. \nonumber 
		\end{align}

  It remains to show the moment conditions on $a_\omega$, which follow immediately from $(M_3)$, since by the choice of $\nu_{\alpha,\beta}$ and $\phi_{\alpha,\beta}^{0}$, the moments of $a_\omega$ are the same as those of $M$. Indeed, for $|\gamma|+|\eta| \leq N_p$,
		\begin{align*}
			&\int{ a_\omega(z) \, (z-z_B)^{\gamma}\ol{(z-z_B)^{\eta}}\Go(z_B,z)dz} \\
            & =  \sum_{|\alpha|+|\beta|\leq N_p}{\nu_{\alpha,\beta}\ \left(|E_0|^{-1} \int_{E_0}{\phi_{\alpha,\beta}^{0}(z)(z-z_B)^{\gamma}\ol{(z-z_B)^{\eta}}\Go(z_B,z)dz} \right)}\\ 
			 &=  \nu_{\gamma,\eta} =\int_{\C^n}{M(z)(z-z_B)^{\gamma}\ol{(z-z_B)^{\eta}}\Go(z_B,z)dz}.
		\end{align*}
		Thus $a_\omega$ is an approximate $(p,s)$-atom.      
\end{proof}

\section{Equivalence Norm of Finite Atomic Decomposition}\label{sec6}

\begin{lemma}\label{fl1}
    Let $0<p\leq 1\leq q\leq \infty$ with $p<q$.
    Suppose $\supp(a)\subset B(c,r)$ with $0<r<1$ with $\|a\|_{L^q}\leq r^{2n(\frac{1}{q}-\frac{1}{p})}$. If there is a $\theta\in \C^n$ and $D>0$ such that $|\theta-c|<D$ and 
    $$\int_{\C^n} a(w)(w-c)^{\Ga}\ol{(w-c)}^\Gb\Go(\theta,w)dw=0, \quad \forall |\Ga|+|\Gb|\leq N_p,$$
    then we can write $a$ as a finite linear combination of $(p,q)$-atoms with $a=\sum_{j=1}^N \Gl_ja_j$ and
    $$\sum_{j=1}^N|\Gl_j|^{p}\leq C_{n,p,q,D}.$$
    The constant $C_{n,p,q,D}$ is independent of $\theta,c,r$.
\end{lemma}
\begin{proof}
The proof is essentially from \cite{BJ}*{Lemma 3.6}. 

Write $A(w)=\Go(\theta,w)a(w)$, $B(w)=\Go(c,w)a(w)$ and $m_{\Ga,\Gb}(a):=\int a(w)(w-c)^{\Ga}\ol{(w-c)}^{\Gb}dw$. Then, $\int A(w)(w-c)^{\Ga}\ol{(w-c)}^{\Gb}dw=0$ for $|\Ga|+|\Gb|\leq N_p$.
Notice that 
$$m_{\Ga,\Gb}(B)=\int A(w)\Go(c-\Gth,w)(w-c)^{\Ga}\ol{(w-c)}^{\Gb}dw.$$
Then, by change of variables, 
$$m_{\Ga,\Gb}(B)=\int A(c+w)w^{\Ga}\ol{w}^{\Gb}\Go(c-\Gth, c+w)dw$$
and using a standard Taylor expansion of $\Go(c-\Gth, c+w)$ with respect to $w$ up to order $N_p-|\Ga|-|\Gb|$, then using the fact that $|\Gth-c|\leq D$ (so that all derivatives are uniformly bounded), the integrals of all terms are 0 except the error term; and thus,
$$|m_{\Ga,\Gb}(B)|\leq C\int_{|w|\leq r}|A(c+w)||w|^{|\Ga|+|\Gb|}|w|^{N_p-|\Ga|-|\Gb|+1}dw\leq C''r^{N_p+1-\Gg_p}.$$

Note that $0<N_p+1-\Gg_p\leq 1$. 

We now proceed as in Lemma \ref{mole1}. Define 
$$\Pi_Rg(z):=\sum_{|\Ga|+|\Gb|\leq N_p}m_{\Ga,\Gb}(g)\frac{1}{R^{|\Ga|+|\Gb|+2n}}\phi_{\Ga,\Gb}\bigg(\frac{z-c}{R}\bigg)=:\sum_{|\Ga|+|\Gb|\leq N_p}m_{\Ga,\Gb}(g)\phi^R_{\Ga,\Gb}(z),$$
where $\phi_{\Ga,\Gb}$ with $\supp(\phi_{\Ga,\Gb})\subset B(0,1)$ satisfies $\int_{B(0,1)} z^{\Ga'}\ol{z}^{\Gb'}\phi_{\Ga,\Gb}(z)dz=\delta_{(\Ga,\Gb)(\Ga',\Gb')}$ and $\delta_{(\Ga,\Gb)(\Ga',\Gb')}$ denotes the Kronecker  delta of the multiindices $\Ga,\Gb$, and $\phi_{\Ga,\Gb}^R\in C^{\infty}_c(\C^n)$.
% For simplicity, we write 
% $$\phi^R_{\Ga,\Gb}:=\frac{1}{R^{|\Ga|+|\Gb|+2n}}\phi_{\Ga,\Gb}\bigg(\frac{z-c}{R}\bigg).$$

Now we write 
$$B=[B-\Pi_{\Ga,\Gb}^{r}B]-\sum_{j=1}^{J}[\Pi_{\Ga,\Gb}^{2^{j}r}B-\Pi^{2^{j-1}r}_{\Ga,\Gb}B]+\Pi^{2^Jr}_{\Ga,\Gb}B=:\Gl_0a'_0+\sum_{j=1}^{J}\Gl_ja'_j+\Gl_{J+1}a'_{J+1},$$
where $J=\inf\{j:2^jr\geq1\}$. The desired decomposition is $a_j=a_j'\Go(w,c)$.

It is easy to see that 
$$\supp([B-\Pi_{\Ga,\Gb}^{r}B])\subset B(c,r), \quad \supp([\Pi_{\Ga,\Gb}^{2^{j}r}B-\Pi^{2^{j-1}r}_{\Ga,\Gb}B])\subset B(c, 2^jr), \quad \supp(\Pi^{2^Jr}_{\Ga,\Gb}B)\subset B(c,2^Jr).$$

First we verify the cancellation for terms $j=0,\cdots,J$. This follows from the construction that
\begin{align*}
    &\int(\Pi_{\Ga,\Gb}^{2^{j}r}B-\Pi^{2^{j-1}r}_{\Ga,\Gb}B)\Go(w,c)[(w-c)^{\Ga'}\ol{(w-c)}^{\Gb'}]\Go(c,w)dw\\
    &=\sum_{\Ga,\Gb}m_{\Ga,\Gb}(B)\int (w-c)^{\Ga'}\ol{(w-c)}^{\Gb'} [\phi^{2^jr}_{\Ga,\Gb}(w)]dw-\sum_{\Ga,\Gb}m_{\Ga,\Gb}(B)\int (w-c)^{\Ga'}\ol{(w-c)}^{\Gb'} [\phi^{2^{j-1}r}_{\Ga,\Gb}(w)]dw\\
    &=m_{\Ga',\Gb'}(B)-m_{\Ga',\Gb'}(B)=0,
\end{align*}
and the cancellation of the term $[B-\Pi_{\Ga,\Gb}^{r}B](w)\Go(w,c)$ can be seen from the definition of projection.

Note that 
$$\|\phi^R_{\Ga,\Gb}(z)\|_{L^q}\lesssim R^{-2n(1-\frac{1}{q})-|\Ga|-|\Gb|}.$$
Thus,
\begin{align*}
\|[\Pi_{\Ga,\Gb}^{2^{j}r}g-\Pi^{2^{j-1}r}_{\Ga,\Gb}g]\|_{L^q}&\leq \sum_{|\Ga|+|\Gb|\leq N_p}C' r^{N_p+1-\Gg_p}(2^{j}r)^{-2n(1-\frac{1}{q})-|\Ga|-|\Gb|}\\
&\leq (2^jr)^{2n(\frac{1}{q}-\frac{1}{p})} \sum_{|\Ga|+|\Gb|\leq N_p}C'r^{N_p+1-|\Ga|-|\Gb|}(2^{j})^{2n(\frac{1}{p}-1)-|\Ga|-|\Gb|};
\end{align*}
and therefore, $|\Gl_j|\leq \sum_{|\Ga|+|\Gb|\leq N_p}C'r^{N_p+1-|\Ga|-|\Gb|}(2^{j})^{2n(\frac{1}{p}-1)-|\Ga|-|\Gb|}$. Thus, the norm can be calculated as
\begin{align*}
    \sum_{j=0}^{J+1}|\Gl_j|^p&\lesssim 1+\sum_{j=1}^{J+1}\sum_{|\Ga|+|\Gb|\leq N_p}r^{(N_p+1-|\Ga|-|\Gb|)p}(2^{j})^{2np(\frac{1}{p}-1)-|\Ga|p-|\Gb|p}\\
    %&\leq C'''\sum_{|\Ga|+|\Gb|\leq N_p}r^{(N_p+1-|\Ga|-|\Gb|)p}(r^{-1})^{2np(\frac{1}{p}-1)-|\Ga|p-|\Gb|p}
    &\simeq 1+r^{p(N_p+1-\Gg_p)}\lesssim 1
\end{align*}
as $|\Ga|+|\Gb|\leq N_p\leq \Gg_p$ whenever $N_p\neq \Gg_p$. If $N_p=\Gg_p$ (\textit{i.e.} $p=\frac{2n}{2n+k}$), then 
\begin{align*}
    \sum_{j=0}^{J+1}|\Gl_j|^p&\lesssim1+ \sum_{j=1}^{J+1}\sum_{|\Ga|+|\Gb|\leq N_p}r^{(N_p+1-|\Ga|-|\Gb|)p}(2^{j})^{2np(\frac{1}{p}-1)-|\Ga|p-|\Gb|p}\\
    %&\leq C'''\sum_{|\Ga|+|\Gb|\leq N_p}r^{(N_p+1-|\Ga|-|\Gb|)p}(r^{-1})^{2np(\frac{1}{p}-1)-|\Ga|p-|\Gb|p}
    &\simeq1+ r^p(1+\log(r^{-1}))\lesssim 1.
\end{align*}
\end{proof}

\begin{lemma}\label{fl2}
    If $F\in \sH^p(\C^n)$ with $\supp(F)\subset B(c,r)$, then the shifting $G_c(z):=F(z)\Go(c,z)\in h^p(\R^{2n})$.
\end{lemma}
\begin{proof}
    Since 
    \begin{align*}
        G_c\ast \phi_s(z)&=\int_{\C^n}F(w)\Go(c,w)\frac{1}{s^{2n}}\phi\bigg(\frac{z-w}{s}\bigg)dw \\
       % &=\int_{\C^n}F(z-tw)\Go(c,z-tw)\phi(w)dw \\
        &=\int_{\C^n}F(z-w)\Go(c,z-w)\phi_s(w)dw\\
        &=\int_{\C^n}F(z-w)[\phi_s(w)\Go(c,z-w)\Go(z,w)]\Go(w,z)dw\\
        &=F\times (\psi_{z,s,c})_s(z).
    \end{align*}
Note that for given $c\in \C^n$, $z\in \C^n$, and $s\in (0,1)$, the function $\psi_{z,s,c}(w):=\phi(w)\Go(c-sw,z-sw)\in C^{\infty}_c(B(0,1))$ because $\phi\in C^{\infty}_c(B(0,1))$. Moreover, $\sup_{|\Ga|\leq N}|\partial^{\Ga}_w\psi_{z,s,c}(w)|\leq C_{\Ga, r}$ using Leibniz rule and the fact that $|z-c|\leq |z-w-c|+|w|\leq r+s<r+1$. Thus, from the hypothesis,
\begin{align*}
    |G_c\ast \phi_s(z)|&\lesssim \sup_{\phi\in \cF_N}\sup_{0<t<1}|F\times \phi_t(z)|\in L^p(\C^n)
\end{align*}
and thus $G_c\in h^p(\R^{2n})$.
% \begin{align*}
%     (F\times\phi_t)(z)&=
% \int G_c(w)\phi_t(z-w)\Go(z-c,w)\,dw\\
% &=\int G_c(z-tw)\phi(w)\Go(z-c, z-tw)dw\\
% &=\Go(z-c,z)\int G_c(z-w)t^{-2n}\phi(w/t)[\Go(c-z,w)]dw
% \end{align*}
%By reversing the argument we have $$\sup_{\phi\in \cF}\sup_{0<t<1}|F\times \phi_t(z)|\simeq \sup_{\phi\in \cF}\sup_{0<t<1}|G_c\ast \phi(z)|.$$
\end{proof}

We now prove the main theorem of this section.

\begin{theorem}\label{ft1}
Let $0<p\leq 1\leq q< \infty$ with $p<q$. If $F\in \sH^{p,q}_{fin}$, then $\|F\|_{\sH^{p,q}_{fin}}\simeq \|F\|_{\sH^p}.$ %On the other hand, if $F\in L^q\cap \mathscr{H}^p(\C^n)$ for some $1<q<\infty$ with support in $B(c,r)$, then $F\in \sH^p_{fin}$.
\end{theorem}
\begin{proof}
    Note that $\|F\|_{\sH^p}\leq \|F\|_{\sH^{p,q}_{fin}}$ as we are limiting the way of doing decompositions. It remains to show that $\|F\|_{\sH^{p,q}_{fin}}\leq C\|F\|_{\sH^p}$.

    Suppose $F$ is a finite combination of  $(p,q)$-atoms. We first consider $q>1$. Take a smooth partition of unity with $\supp(\Gz_j)\subset B(z_j,1)$ and $\sum_j\Gz_j=1$. There exists a finite set $J_0$ such that  $F=\sum_{j\in J_0} \Gz_jF$. Now we write $G_j(w)=F(w)\Gz_j(w)\Go(z_j,w)$. Then, 
    $$\sum_{j\in J_0}\|G_j\|_{h^p}^p\leq C\|f\|_{\sH^p}^p.$$  
    
    Since $G_j\in h^p\cap L^q$ with $\supp(G_j)\subset B(z_j,1)$, by \cite{YY}*{Theorem 6.3 (i) and proof of Case 1} we obtain a finite atomic decomposition such that
    $$G_j=\sum_{\nu\in S_j}\mu_{j,\nu}a_{j,\nu},$$
    where $a_{j,\nu}$ are classical $(p,q)$-atoms supported in $B(z_j,K)$ for fixed $K>1$, and  
    $$\sum_{\nu\in S_j}|\mu_{j,\nu}|^p\leq C\|G_j\|_{h^p}^p.$$

    Now note that  
    $$F(z)\Gz_j(z)=\sum_{\nu\in S_j}\mu_{j,\nu}a_{j,\nu}(z)\Go(z,z_j),$$
    Note that $$\int_{\C^n} a_{j,\nu}(w)\Go(w,z_j)(w-c_{j,\nu})^{\Ga}\ol{(w-c_{j,\nu})}^{\Gb}\Go(z_j,w)dw=0.$$
    We apply Lemma \ref{fl1} to atoms with $r_{j,\nu}<1$ to obtain the further decomposition of $a_{j,\nu}(w)\Go(w,z_j)$, and define $M_{j,\nu}=1, \ \Gg_{j,\nu,1}=1, \ a_{j,\nu,m}=a_{j,\nu}\Go(\cdot,z_j)$ for $r_{j,\nu}\geq 1$, then we have 
    $$a_{j,\nu}(w)\Go(w,z_j)=\sum_{m=1}^{M_{j,\nu}}\Gg_{j,\nu,m}a_{j,\nu,m}\text{ with }  \sum_{m=1}^{M_{j,\nu}}|\Gg_{j,\nu,m}|^p\leq C.$$

    Thus, we have a finite decomposition of $F$
    $$F=\sum_{j\in J_0}\sum_{\nu\in S_j}\sum_{m=1}^{M_{j,\nu}}\mu_{j,\nu}\Gg_{j,\nu,m}a_{j,\nu,m},$$
    using $(p,q)$ atoms, and 
    \begin{align*}
        \|F\|_{\sH^{p,q}_{fin}}^p\leq\sum_{j\in J_0}\sum_{\nu\in S_j}\sum_{m=1}^{M_{j,\nu}} |\mu_{j,\nu}\Gg_{j,\nu,m}|^p\leq C \sum_{j\in J_0}\sum_{\nu\in S_j}|\mu_{j,\nu}|^p\leq C'\|F\|_{\sH^p}^p.
    \end{align*}
    This completes the proof for $q>1$.

    We now consider the case $q=1$ and $0<p<1$. Let $F\in \sH^{p,1}_{fin}$ with $\supp(F)\subset B(c,r)$. Then, for $j\in \N$, we write $F_j=F\chi_{\{|F|\leq j\}}$ and $R_j=F-F_j$. We can see that $R_j$ is a constant multiple of a $(p,1)$-atom because $\|R_j\|_{L^1}=\|F-F_j\|_{L^1}\leq 2\|F\|_{L^1}$, and we enlarge its supporting ball at least 1. Note that $\|R_j\|_{L^1}\to0$ as $j\to\infty$ by Dominated Convergence Theorem.  On the other hand, since $F_j\in L^{\infty}$, then $F_j$ is a constant multiple of a $(p,2)$-atom, and therefore, by the previous case, we have $\|F_j\|_{\sH^{p,2}_{fin}}\simeq \|F_j\|_{\sH^p}$. Let $\Ge>0$. If we choose $j$ (depending on $r$) sufficiently large so that $\|R_j\|_{\sH^{p,1}_{fin}}\leq \Ge$, we have 
$$\|F\|_{\sH^{p,1}_{fin}}^p\leq \|F_j\|_{\sH^{p,1}_{fin}}^p+\|R_j\|_{\sH^{p,1}_{fin}}^p\lesssim \|F_j\|_{\sH^p}^p+\Ge^p\lesssim \|F\|_{\sH^p}^p+2\Ge^p.$$
As $\Ge$ is arbitrary, $\|F\|_{\sH^{p,1}_{fin}}\lesssim\|F\|_{\sH^p}$.
\end{proof}

\begin{remark}
    If we use $(p,q)$-atoms with extra cancellation up to $|\Ga|+|\Gb|\leq N$ for any integer $N\geq N_p$, we still have the equivalence between atomic and finite atomic quasinorms for $p<q<\infty$.

    Indeed, applying the projection construction in the proof of Lemma \ref{fl1} with order $N$ to a $(p,q)$-atom, we obtain
    $$
    \sum_{j=0}^{J+1}|\Gl_j|^p
    \lesssim 1+
    \sum_{k=N_p+1}^{N}\sum_{j=0}^{J}
    (2^j)^{p(\Gg_p-k)}
    \lesssim 1,
    $$
    where we have used the fact that the moments of order $k>N_p$ are bounded by $Cr^{k-\Gg_p}$ and that $N_p+1>\Gg_p$.
\end{remark}

Next, if we consider a finite approximate atomic decomposition, do we have the norm equivalence?

\begin{lemma}\label{fl3}
    Suppose $a$ is an approximate $(p,q)$-atom. Then, there exists a finite sequence of $(p,q)$-atom $\{a_j\}_{j=0}^{J+1}$ and a finite sequence of scalar $\{\Gl_j\}_{j=0}^{J+1}$ such that 
    $$a=\sum_{j=0}^{J+1}\Gl_ja_j$$
    and $\sum_{j=0}^{J+1}|\Gl_j|^p\lesssim_p 1,$
    where the implicit constant is independent of $a$.
\end{lemma}
\begin{proof}
Let $\supp(a)\subset B(c,r)$. If $r\geq 1$, it is immediate by taking itself as the decomposition, as no moment condition is needed. Following the notation in Lemma \ref{fl1}, note that if $B(z)=a(z)\Go(c,z)$, then 
    $$B=(B-\Pi_{\Ga,\Gb}^rB)-\sum_{j=1}^{J}[\Pi^{2^jr}_{\Ga,\Gb}B-\Pi^{2^{j-1}r}_{\Ga,\Gb}B]+\Pi^{2^Jr}_{\Ga,\Gb}B=:\Gl_0a_0'+\sum_{j=1}^{J}\Gl_ja_j'+\Gl_{J+1}a'_{J+1}.$$
The desired decomposition is obtained by setting  $a_j(z)=a'_j(z)\Go(z,c)$ once we verify the $p$-sum of the coefficients.

If $\Gg_p\neq N_p$, then 
$|\Gl_0|\lesssim 1$, $|\Gl_{j}|\lesssim C_p(2^{j}r)^{\Gg_p-N_p}$ for $j\in \{1,\cdots, J\}$, and $|\Gl_{J+1}|\lesssim_p1$ using $2^{J-1}r\leq 1\leq 2^{J}r$.
Then, 
$$\sum_{j=0}^{J+1}|\Gl_j|^p\lesssim1+r^{p(\Gg_p-N_p)}\sum_{j=0}^{J}(2^{j})^{p(\Gg_p-N_p)} \lesssim 1.$$

If $\Gg_p=N_p$, then 
\begin{align*}
    |\Gl_j|&\lesssim \sum_{|\Ga|+|\Gb|\leq N_p-1} (2^{j}r)^{2n(\frac{1}{p}-1)-|\Ga|-|\Gb|}+\sum_{|\Ga|+|\Gb|=N_p}\frac{1}{[\log(1+r^{-1})]^{\frac{1}{p}}}\\
    &\lesssim 2^jr+ \frac{1}{[\log(1+r^{-1})]^{\frac{1}{p}}},
\end{align*}
while $|\Gl_0|,|\Gl_{J+1}|\lesssim_p 1$; and thus,
$$\sum_{j=0}^{J+1}|\Gl_j|^p\lesssim 2+r^p \sum_{j=1}^{J}2^{jp}+\frac{J}{\log(1+r^{-1})}\lesssim 2+C_p(2^{J-1}r)^p+C_p'\frac{\log(r^{-1})}{\log(1+r^{-1})}\lesssim_p1.$$
\end{proof}

\begin{theorem}
    Let $0<p\leq1\leq q<\infty$ with $p<q$. If $f$ can be written as a finite linear combination of approximate $(p,q)$-atoms, then $f\in \sH^{p,q}_{fin}$, and thus 
    $$\|f\|_{\sH^p}\simeq \|f\|_{\sH^{p,q}_{fin}}\simeq \inf\bigg(\sum_{j=1}^{N}|\Gl_j|^p\bigg)^{\frac{1}{p}},$$
    where the infimum is taken over all possible finite linear combinations of approximate $(p,q)$-atoms.
\end{theorem}
\begin{proof}
    For convenience, write $\|f\|_{\sH^{p,q}_{ap,fin}}=\inf\bigg(\sum_{j=1}^{N}|\Gl_j|^p\bigg)^{\frac{1}{p}}$, where the infimum is taken over all possible finite linear combinations of approximate $(p,q)$-atoms.

    Then, if one writes $f=\sum_{k=1}^{K}\Gl_ka_k$, where $a_k$ are approximate $(p,q)$-atoms, then by Lemma \ref{fl3},
    $$f=\sum_{k=1}^{K}\sum_{j=0}^{J_k+1}\Gl_k\Gl_k^ja_{j,k},$$
    where $a_{j,k}$ are $(p,q)$-atoms; and thus,
    $$\|f\|_{\sH^{p,q}_{fin}}^p\lesssim \sum_{k=1}^{K}\sum_{j=0}^{J_{k}+1}|\Gl_k\Gl_k^j|^p\lesssim \sum_{k=1}^K|\Gl_k|^p.$$
    Taking the infimum on both sides over finite linear combinations of approximate $(p,q)$-atoms, we have 
    $$\|f\|_{\sH^{p,q}_{fin}}\lesssim \|f\|_{\sH^{p,q}_{ap, fin}}.$$
    On the other hand, since every atom is an approximate atom, we have 
    $$\|f\|_{\sH^{p,q}_{ap,fin}}\lesssim \|f\|_{\sH^{p,q}_{fin}}.$$
    Thus, $\|f\|_{\sH^{p,q}_{ap, fin}}\simeq \|f\|_{\sH^{p,q}_{fin}}\simeq \|f\|_{\sH^p}.$
\end{proof}
We remark that this result is also true for $h^p(\R^n)$. 

\section{Singular integral operators on twisted Hardy spaces}\label{sec4}
\subsection{Singular integral operators}
We first define the singular integral operators with respect to twisted convolution. If $K\in \mathcal{D}'(\C^n\setminus\{0\})$, for $f\in C^{\infty}_c(\C^n)$, we define
$$T_Kf(z):=f\times K(z)$$
for suitable assumptions on $K$.

Here we impose the following assumptions:
\begin{enumerate}
    \item the $L^2$ boundedness of $K$: 
    \begin{align}\label{L2K}
    \|f\times K\|_{L^2(\C^n)}\leq C_K\|f\|_{L^2(\C^n)},
\end{align}
\item the size condition on $K$:
\begin{align} \label{sizeK}
    |K(z)|\leq C_{K}'\min\{|z|^{-2n}, |z|^{-2n-\eta}\}
\end{align}
\item the smoothness condition on $K$: there exist
$s\in[1,2]$, $N_0\in\N\cup\{0\}$ and
$\Gd\in(N_0,N_0+1]$ such that for every
$r\in(0,1)$ and $j\in\N$, there exists a
polynomial $P_z\in\sP_{N_0}$, whose coefficients
are measurable with respect to $z$, satisfying
\begin{align}\label{smoothK}
    \int_{2^jr\leq |z|\leq 2^{j+1}r}
    |K(z-w)\Go(z,w)-P_z(w)|^s dz
    \leq C'(2^j)^{2n(1-s)-s\Gd}
    r^{2n(1-s)}
\end{align}
uniformly for all $w\in B(0,r)$;
\item cancellation condition: for all $B\subset \C^n$ with $r(B)<1$, for all multiindices $|\Ga|+|\Gb|\leq N_p$, the element $f(z):=T_K^*((z-c_B)^{\Ga}\ol{(z-c_B)}^{\Gb}\Go(z,c_B))$ satisfies the following: there exists a polynomial $P^{N_p}_B(f)(z)\in \sP_{N_p}$ such that 
\begin{align}\label{can1}
    \bigg(\frac{1}{|B|}\int_B|f(z)\Go(c_B,z)-P^{N_p}_B(f)(z)|^2dz\bigg)^{\frac{1}{2}}\leq C'[r(B)]^{\Gg_p}\phi_{p,\Ga,\Gb}(r(B)),
\end{align}
where $\phi_{p,\Ga,\Gb}(t)=1$ if $|\Ga|+|\Gb|<\Gg_p$ and  $\phi_{p,\Ga,\Gb}(t)=[\log(1+1/r(B))]^{-\frac{1}{p}}$ if $|\Ga|+|\Gb|=\Gg_p=N_p$.
\end{enumerate}

We remark that \eqref{smoothK} can also be written in the form centered at $c$:
$$ \int_{2^jr\leq |z-c|\leq 2^{j+1}r}|K(z-w)\Go(z-c,w-c)-P_{z-c}(w-c)|^s dz \leq C' (2^j)^{2n(1-s)-s\delta}r^{2n(1-s)}$$
for all $w\in B(c,r)$ with $0<r<1$.

% \begin{remark}
%     We first remark that with the assumption \eqref{sizeK}, \eqref{smoothK} is indeed equivalent to \begin{align}\label{smoothK1}
%     \int_{2^jr\leq |z|\leq 2^{j+1}r}|K(z-w)-P_z(w)|^s dz \leq C'' (2^j)^{2n(1-s)-s\delta}r^{2n(1-s)}.
% \end{align}
% Indeed, 
% \begin{align*}
%     &\int_{2^jr\leq |z-c|\leq 2^{j+1}r}|K(z-w)|^s|\Go(z-w)-1|^sdz \\
%     &\lesssim  \min\{(2^{j}r)^{-2ns}, (2^jr)^{-2ns-\eta s}\}(2^{j}r)^sr^s(2^{j+1}r)^{2n}\\
%     &\lesssim \min\{(2^j)^{2n(1-s)+s}r^{2s+2n(1-s)}, (2^j)^{2n(1-s)+s-\eta s}r^{2s+2n(1-s)-\eta s}\}
% \end{align*}
% \end{remark}

% \begin{lemma}
%     Let $T_K$ defined above satisfy the assumptions \eqref{L2K}-\eqref{smoothK} with $s=1$. Then, $(z-c)^{\Ga}\ol{(z-c)}^{\Gb}\cdot\allowbreak T_Kg(z)\in L^1(\C^n)$ for all $g\in L^2(\C^n)$ that there exists $c\in \C^n$, $r>0$ such that $\supp(g)\subset B(c,r)$ and $\int g(w)(w-c)^{\Ga}\ol{(w-c)}^{\Gb}\Go(c,w)dw=0$ for all $|\Ga|+|\Gb|\leq N$ provided that $N<\min\{\eta,\delta\}$.
% \end{lemma}
\begin{lemma}
    Let $T_K$ satisfy \eqref{L2K} and \eqref{sizeK}. Suppose that $g\in L^2(\C^n)$ with $\supp(g)\subset B(c,r)$ for some $c\in\C^n$ and $r>0$. Then $$(z-c)^\Ga\ol{(z-c)}^\Gb T_Kg(z)\in L^1(\C^n)$$
    provided that $|\Ga|+|\Gb|<\eta$.
\end{lemma}
\begin{proof}
    Let $g$ be given in the statement. Write $R=\max\{1,r\}$,
    $$\int_{\C^n}|(z-c)^{\Ga}\ol{(z-c)}^{\Gb}T_Kg(z)|dz=\int_{B(c,2R)}+\int_{[B(c,2R)]^c}|(z-c)^{\Ga}\ol{(z-c)}^{\Gb}T_Kg(z)|dz.$$

    From the boundedness of $T_K$ on $L^2(\C^n)$ we have 
    $$\int_{B(c,2R)}|(z-c)^{\Ga}\ol{(z-c)}^{\Gb}T_Kg(z)|dz\leq R^{|\Ga|+|\Gb|+n}\|g\|_{L^2}<\infty.$$

    To estimate the second term, first note that 
    $$|z-w|\geq |z-c|-r\geq \frac{1}{2}|z-c|\geq 1$$
provided that $w\in B(c,r)$ and $z\notin B(c,2R)$.
Using \eqref{sizeK}, we have 
$$|T_Kg(z)|\leq \int_{B(c,r)}|g(w)||K(z-w)|dw\lesssim |z-c|^{-2n-\eta}\|g\|_{L^1}.$$
Therefore, 
\begin{align*}
    \int_{[B(c,2R)]^c}|(z-c)^{\Ga}\ol{(z-c)}^{\Gb}T_Kg(z)|dz&\lesssim \|g\|_{L^1}\int_{|z-c|\geq R} |z-c|^{k-2n-\eta}dz \\
    &\lesssim \frac{1}{\eta-k}R^{k-\eta}\|g\|_{L^1}<\infty
\end{align*}
provided that $k=|\Ga|+|\Gb|<\eta$. 
%     If $r\geq 1$, then  
%     \begin{align*}
%         &\int_{B(c,2r)}|(z-c)^{\Ga}\ol{(z-c)}^{\Gb}T_Kg(w)|dw\\
%         &\leq \sum_{j=1}^{\infty}\int_{2^jr\leq |w-c|\leq 2^{j+1}r}\int_{B(c,r)}|K(w-z)||g(z)| (2^{j+1}r)^{|\Ga|+|\Gb|} dzdw \\
%         &\lesssim_{\Ga,\Gb} r^{n+|\Ga|+|\Gb|}\|g\|_{L^2}\sum_{j=1}^{\infty}(2^j)^{|\Ga|+|\Gb|}\int_{2^jr\leq |w-c|\leq 2^{j+1}r}|w-c|^{-2n-\eta}dz\\
%          &\lesssim_{\Ga,\Gb} r^{n+|\Ga|+|\Gb|-\eta}\|g\|_{L^2}\sum_{j=1}^{\infty}(2^j)^{|\Ga|+|\Gb|-\eta}<\infty
%     \end{align*}
%     provided that $|\Ga|+|\Gb|\leq N<\eta$.

%     If $r<1$, note that \begin{align*}
%     T_Ka(z)&=\int_{\C^n}a(w)K(z-w)\Go(z,w)dw\\
%     &=\int_{\C^n} a(w)[K(z-w)\Go(z,w)\Go(w,c)-P(w-c)]\Go(c,w)dw\\
%     &= \Go(z,c)\int_{B(c,r)}a(w)[K(z-w)\Go(z-c,w-c)-P(w-c)]\Go(c,w)dw;
% \end{align*}
%  then the cancellation of $g$ gives 
%         \begin{align*}
%         &\int_{B(c,2r)}|(w-c)^{\Ga}\ol{(w-c)}^{\Gb}T_Kg(w)|dw\\
%         &\leq \sum_{j=1}^{\infty}\int_{2^jr\leq |w-c|\leq 2^{j+1}r}\int_{B(c,r)}|K(w-z)\Go(w-c,z-c)-P(z-c)||g(z)| (2^{j+1}r)^{|\Ga|+|\Gb|} dzdw \\
%         &\lesssim r^{n+|\Ga|+|\Gb|}\|g\|_{L^2}\sum_{j=1}^{\infty}(2^j)^{|\Ga|+|\Gb|-\Gd}<\infty
%     \end{align*}
%     provided that $|\Ga|+|\Gb|\leq N<\Gd$.
\end{proof}

As a consequence of the preceding lemma, for every twisted polynomial $Q$ of degree strictly less than $\eta$, the mapping
$$a\mapsto \langle T_Ka,Q\rangle$$
defines a bounded linear functional on $L^2(B)$ for every ball $B\subset\C^n$. Hence, by the Riesz representation theorem, $T_K^*Q$ is well-defined as an element of $L^2_{\mathrm{loc}}(\C^n)$ satisfying $$
\langle T_Ka,Q\rangle=\langle a,T_K^*Q\rangle
$$
for all $a\in L_c^2(\C^n)$. This justifies the notation in \eqref{can1}.

\begin{theorem}\label{main1}
   Under the assumptions of \eqref{L2K}-\eqref{can1}, the operator $T_K$ maps $(p,2)$ atoms satisfying
(A3) up to order
$N_0$ to approximate $(p,s,\Gl)$ molecule provided that $p<s$ and 
$$\Gg_p<\frac{\Gl}{s}-\frac{2n}{s'}<\min\{\eta,\Gd\}.$$
Consequently, $T_K$ can be extended to a bounded operator on $\sH^p(\C^n)$.
\end{theorem}
\begin{proof}
By Remark \ref{rem:HigherCancellation} and
Theorem \ref{ft1}, it suffices to consider a $(p,2)$-atom $a$ supported in $B(c,r)$ satisfying (A3) up to order $N_0$.

    It suffices to show that $T_K(a)$ is an $(p,s,\Gl)$ for a $(p,2)$ atom $a$ with support in $B(c,r)$.

    We first show (M1) regardless of $a$. Note that 
    \begin{align*}
        \|T_Ka\|_{L^s(2B)}&\leq \|T_Ka\|_{L^2(2B)}|2B|^{\frac{1}{s}-\frac{1}{2}}\lesssim \|K\|r^{2n(\frac{1}{s}-\frac{1}{p})}.
    \end{align*}

    We now show (M2). If $r\geq 1$, then,
    \begin{align*}
    \|T_Ka|\cdot -c|^{\frac{\Gl}{s}}\|_{L^s([B(c,2r)]^c)}^s&=\int_{|z-c|\geq 2r} |z-c|^{\Gl} |a\times K(z)|^s dz \\
    &\lesssim \int_{|z-c|\geq 2r}|z-c|^{\Gl}r^{-s\Gg_p}|z-c|^{-2sn-s\eta}dz \\
    &\lesssim r^{-s\Gg_p}\int_{|z-c|\geq 2r}\frac{1}{|z-c|^{2sn+s\eta-\Gl}}dz \leq C'' r^{-s\Gg_p}r^{-2sn-s\eta+\Gl+2n}\\
    &\lesssim r^{\Gl+2n(1-\frac{s}{p})-s\eta}\lesssim r^{\Gl+2n(1-\frac{s}{p})}.
\end{align*}
    provided that $-2sn-s\eta+\Gl+2n<0$, that is, $\eta>\frac{\Gl}{s}-\frac{2n}{s'}$.

    If $r<1$, then
\begin{align*}
    T_Ka(z+c)&=\int_{\C^n}a(w)K(z+c-w)\Go(z+c,w)dw\\
    &=\int_{\C^n} a(w)[K(z+c-w)\Go(z+c,w)\Go(w,c)-\Go(z,c)P_z(w-c)]\Go(c,w)dw\\
    &= \Go(z,c)\int_{B(0,r)}a(w+c)[K(z-w)\Go(z,w)-P_z(w)]\Go(c,w)dw
\end{align*}
Thus, if $z\in [B(c,2r)]^c$, then $z-c\in [B(0,2r)]^c$ and
$$T_Ka(z+c)=\Go(z,c)\int_{B(0,r)}a(w+c)[K(z-w)\Go(z,w)-P_z(w)]\Go(c,w)dw.$$
Therefore,
    \begin{align*}
    &\|T_Ka|\cdot -c|^{\frac{\Gl}{s}}\|_{L^s([B(c,2r)]^c)}^s\\
    &=\|T_Ka(\cdot+c) |\cdot|^{\frac{\Gl}{s}}\|_{L^s([B(0,2r)]^c)}^s\\
    &=\int_{|w|\geq 2r} |z|^{\Gl} |T_Ka(z+c)|^s dz \\
    &\leq \sum_{j=1}^{\infty}\int_{2^{j}r\leq |z|\leq 2^{j+1}r} \bigg(\int_{B(0,r)}|K(z-w)\Go(z,w)-P_z(w)||a(w+c)||z|^{\frac{\Gl}{s}}dw\bigg)^sdz\\
    &\leq \sum_{j=1}^{\infty}(2^{j+1}r)^{\Gl}\bigg(\int_{B(0,r)}|a(w+c)|\bigg[\int_{2^{j}r\leq |z|\leq 2^{j+1}r}|K(z-w)\Go(z,w)-P_z(w)|^s dz\bigg]^{\frac{1}{s}}dw\bigg)^s\\
    &\leq \sum_{j=1}^{\infty}(2^{j+1}r)^{\Gl}(2^jr)^{-2n(s-1)}2^{-js\Gd}r^{-s\Gg_p}\lesssim r^{\Gl+2n(1-\frac{s}{p})}
    \end{align*} 
    % \begin{align*}
    % &\|a\times K|\cdot -c|^{\frac{\Gl}{s}}\|_{L^s([B(c,2r)]^c)}^s\\
    % &=\int_{|z-c|\geq 2r} |z-c|^{\Gl} |a\times K(z)|^s dz \\
    % &=\sum_{j=0}^{\infty}\int_{2^{j}r\leq |z-c|\leq 2^{j+1}r} \bigg(\int_{B(c,r)}|K(z-w)\Go(z,w)-P_c(z)||a(w)||z-c|^{\frac{\Gl}{s}}dw\bigg)^sdz\\
    % &\leq \sum_{j=0}^{\infty}(2^{j+1}r)^{\Gl}\bigg(\int_{B(c,r)}|a(w)|\bigg[\int_{2^{j}r\leq |z-c|\leq 2^{j+1}r}|K(z-w)\Go(z,w)-P_c(z)|^s dz\bigg]^{\frac{1}{s}}dw\bigg)^s\\
    % &\leq \sum_{j=0}^{\infty}(2^{j+1}r)^{\Gl}(2^jr)^{-2n(s-1)}2^{-js\Gd}r^{-s\Gg_p}\lesssim r^{\Gl+2n(1-\frac{s}{p})}
    % \end{align*}
    provided that $\Gl<2n(s-1)+s\Gd$.

    We now show (M3'). Note that if one writes $f=T^*_K((\cdot-c)^{\Ga}\ol{(\cdot-c)^{\Gb}}\Go(\cdot,c))$, 
    \begin{align*}
        \bigg|\int_{\C^n} T_Ka(z)(z-c)^{\Gb}\ol{(z-c)}^{\Ga}\Go(c,z)dz\bigg|&=|\langle T_Ka, \ol{(\cdot-c)^{\Gb}}(\cdot -c)^{\Ga}\Go(\cdot,c)\ran|\\
        &\leq r^{2n(\frac{1}{2}-\frac{1}{p})}\bigg(\int_{B(c,r)}|f(z)\Go(c,z)-P_B^{N_p}(f)(z)|^2 dz\bigg)^{\frac{1}{2}}\\
        &\leq C' \phi_{p,\Ga,\Gb}(r),
    \end{align*}
    which is enough to conclude (M3').

    Therefore, we have $\|Ta\|_{\sH^p}\lesssim 1$ and by equivalence of finite atomic norm (see Theorem \ref{ft1}) and the density of finite atomic space, we can conclude that there is a bounded extension of $T_K$ such that $T_K:\sH^p(\C^n)\to\sH^p(\C^n)$. 
\end{proof}

\begin{example}
    In the paper \cite{MPR}, the twisted Riesz transform has the kernel 
$$K(z)=\text{p.v. }\frac{z_j}{|z|^{2n+1}}\phi(z),$$
where $\phi$ is real and radial with $\supp(\phi)\subset B(0,2)$ and $\phi\equiv 1$ on $B(0,1)$. 

It is proved in \cite{MPR} that $T_K$ is bounded on $L^2(\C^n)$. Also, since the kernel has compact support, \eqref{sizeK} is satisfied. We now show \eqref{smoothK}, and indeed, for $4\geq |z|\geq 2r\geq |w|$
\begin{align*}
|K(z-w)\Go(z,w)-K(z)|&\leq |K(z-w)-K(z)||\Go(z,w)|+|K(z)||\Go(z,w)-1|\\
&\leq C_{\phi}\frac{|w|}{|z|^{2n+1}}+\frac{1}{|z|^{2n}}|w||z|\\
&\lesssim \frac{|w|}{|z|^{2n+1}}
\end{align*}
thanks to the compact support of $\phi$. Thus, a standard argument shows that \eqref{smoothK} holds with $P_z(w)=K(z)$ and $\Gd=1$.

Next, we can see that 
\begin{align*}
    T^*_K(\Go(z,c_B))(w)=\Go(w,c_B) \text{ p.v.}\int_{\C^n} \frac{(\ol{u})_j}{|u|^{2n+1}}\phi(u)\Go(u,c_B-w)du=:\Go(w,c_B)M(w-c_B),
\end{align*}
which is impossible to be 0; and since $M(0)=0$, 
$$|M(w-c_B)|\lesssim \int_{B(0,2)}\frac{|\Go(u,c_B-w)-\Go(0,c_B-w)|}{|u|^{2n}}dz\lesssim \int_{B(0,2)}\frac{|c_B-w|}{|u|^{2n-1}}du\simeq |c_B-w|.$$
Thus, 
$$\bigg(\frac{1}{r^{2n}}\int_{B(c_B,r)}| T^*_K(\Go(z,c_B))(w)|^2dw\bigg)^{\frac{1}{2}}\leq C' r\leq C'[\log(1+r^{-1})]^{-1}.$$ 

Thus, $T_K$ is bounded on $\sH^p(\C^n)$ for $\frac{2n}{2n+1}<p\leq 1$, which proves the same result as \cite{HL}*{Theorem 9}.
\end{example}

Note that the theory can be extended to singular integral operators of non-convolution type by defining 
$$T_Kf(z):=\int_{\C^n}K(z,w)f(w)dw$$
for $f\in C^{\infty}_c(\C^n)$, and we replace \eqref{sizeK} by 
\begin{align} \label{sizeK'}
    |K(z,w)|\leq C'_K\min\{|z-w|^{-2n}, |z-w|^{-2n-\eta}\} \tag{4.2'}
\end{align}
for $z\neq w$, and \eqref{smoothK} by  
\begin{align}\label{smoothK'}
    \int_{2^jr\leq |w-c|\leq 2^{j+1}r}|K(w,y)\Go(y,c)-P_{w,c}(y)|^s dw \leq C' (2^j)^{2n(1-s)-s\delta}r^{2n(1-s)} \tag{4.3'}
\end{align} 
with the existence of polynomial $P_{w,c}\in \sP_{N_0}$ and $N_0<\Gd\leq N_0+1$. This is the correct smoothness condition because 
\begin{align*}
    \int K(w,y)a(y)dy&=\int K(w,y)\Go(y,c)a(y)\Go(c,y)dy \\
    &=\int [K(w,y)\Go(y,c)-P_c(y)]a(y)\Go(c,y)dy.
    %&=\int 
    %[K(w,y)-P_c(w)\Go(c,y)]a(y)dy.
\end{align*}

A similar argument as in Theorem \ref{main1} gives the following.

\begin{theorem}
    Suppose that $K(z,w)$ satisfies \eqref{L2K}, \eqref{sizeK'}, \eqref{smoothK'}, and \eqref{can1}.
    If $p<s$ and  $$\Gg_p<\min\{\eta,\Gd\},$$ then the singular integral operator associated with the kernel $K$ can be extended to a bounded linear operator on $\sH^p(\C^n)$.
\end{theorem}

\subsection{Strongly singular integral operators}
In this section, we consider the strongly singular operator of the form 
$$T_Kf(z):=\int K(z,w)f(w)dw,$$
where $K\in \cD'(\C^n\setminus\Delta)$, $\Delta=\{(z,w)\in \C^{2n}:z=w\}$, and $f\in C^{\infty}_c(\C^n)$. 
We impose the following assumptions for fixed $0<\Gs< 1$:
\begin{enumerate}
\item \eqref{sizeK'} and \eqref{can1} hold;
    \item $T$ can be extended to a bounded operator on $L^2(\C^n)$ and there exists $\Gg\in [(1-\Gs)n,n)$ such that $T$ can be also extended to a bounded operator from $L^q(\C^n)$ to $L^2(\C^n)$, where $\frac{1}{q}=\frac{1}{2}+\frac{\Gg}{2n}$; here we write them in short as
    \begin{align}\label{bdd2}
        \|T_Kf\|_{L^2}\leq C_K(\min\{\|f\|_{L^2}, \|f\|_{L^q}\});
    \end{align}
    \item the smoothness condition on $K$: there exists $s\in[1,2]$, $N_0\in \N\cup\{0\}$, and $\Gd\in (N_0, N_0+1]$ such that for $w,c\in \C^n$, there exists $P_{w,c}\in \sP_{N_0}$ with measurable coefficients with respect to $w$ such that 
    \begin{align}\label{smoothK2}
        \bigg(\int_{2^jr^\rho\leq |w-c|\leq 2^{j+1}r^\rho}|K(w,y)\Go(y,c)-P_{w,c}(y)|^s dw\bigg)^{\frac{1}{s}}\leq C' |B(w,2^jr^{\rho})|^{\frac{1}{s}-1+\frac{\Gd}{2n}(\frac{1}{\Gr}-\frac{1}{\Gs})}2^{-\frac{j\Gd}{\rho}}
    \end{align}
    for all $j\in \N$, $r\in (0,1)$, $0<\Gr\leq \Gs$, $c\in \C^n$, $y\in B(c,r)$.
\end{enumerate}

\begin{theorem}
    Suppose $K(z,w)$ satisfies \eqref{bdd2}, \eqref{sizeK'}, \eqref{smoothK2}, and \eqref{can1} with some $0<\Gs<1$, $\Gd>0$, $1\leq s\leq 2$ with $p<s$. Then there exists a bounded extension of $T_K$ on $\sH^p(\C^n)$ provided that
    $$\max\bigg\{\frac{2n}{2n+\eta},p_0\bigg\}<p\leq 1, \quad \text{where }\frac{1}{p_0}=\frac{1}{2}+\frac{\Gg(n+\Gd/\Gs)}{2n(\Gg-\Gd+\Gd/\Gs)}.$$
\end{theorem}
\begin{proof}
The argument is essentially the same as the proof of \cite{DLPV1}*{Theorem 5.8}. By Remark \ref{rem:HigherCancellation}, it suffices to consider a $(p,2)$-atom satisfying the cancellation conditions in (A3) up to order $\max\{N_0,N_p\}$.

Again, we start with a $(p,2)$ atom with support $B=B(c,r)$. We show that $Ta$ is a $(p,s,\Gl)$ approximate molecule for $\Gl$ with $\Gg_p<\frac{\Gl}{s}-\frac{2n}{s'}<\min\{\eta,\Gg_{p_0}\}$, where $\Gg_{p_0}=2n(\frac{1}{p_0}-1)$.

Condition (M1) can be shown as in Theorem \ref{main1} without splitting the cases. For (M2), the same proof works for $r\geq1$. We still need to show the case $r<1$. Instead of considering $\int_{2B}$ and $\int_{(2B)^c}$, we will consider $\int_{2B(c,r^{\rho})}$ and $\int_{[2B(c,r^{\rho})]^c}$. For short, we write $2B^\rho:=B(c,2r^{\rho})$. We consider 
\begin{align*}
    \int_{\C^n}|T_Ka(z)|^s|z-c|^{\Gl}dz &= \int_{2B^\rho}|T_Ka(z)|^s |z-c|^{\Gl}dz+\int_{[2B^{\rho}]^c}|T_Ka(z)|^s |z-c|^{\Gl}dz =:I_1+I_2.
\end{align*}
Using the boundedness of $T$ from $L^q$ to $L^2$, we can estimate $I_1$ as 
\begin{align*}
    I_1&\lesssim r^{\Gl\rho} |2B^{\rho}|^{1-\frac{s}{2}}\|Ta\|_{L^2}^s \lesssim r^{\rho\Gl+2n\rho(1-s/2)}\|a\|_{L^q}^s\lesssim r^{\rho \Gl+2n\rho(1-\frac{s}{2})+2ns(\frac{1}{q}-\frac{1}{p})}\lesssim r^{\Gl+2n(1-\frac{s}{p})}
\end{align*}
provided that $\Gl\leq -2n(1-\frac{s}{2})+\frac{2ns}{1-\rho}(\frac{1}{q}-\frac{1}{2})$. 

To estimate $I_2$, similar to the singular integral operator,  by using the expression of $T_Ka$ and writing splitting the domain into annuli, we have 
\begin{align*}
    I_2&\leq \sum_{j=1}^{\infty}(2^jr^{\rho})^\Gl \bigg\{\int_B|a(y)|\bigg[\int_{2^jr^\rho\leq |w-c|\leq 2^{j+1}r^{\rho}}|K(w,y)\Go(y,c)-P_c(y)|^sdw\bigg]^{\frac{1}{s}}dy\bigg\}^{s}\\
    &\lesssim \sum_{j=1}^{\infty}(2^jr^{\rho})^\Gl \bigg((2^jr^\rho)^{2n(\frac{1}{s}-1)+\Gd(\frac{1}{\rho}-\frac{1}{\Gs})}2^{-j\frac{\Gd}{\rho}}\bigg)^s \|a\|_{L^2}^s|B|^{\frac{s}{2}}\\
    &\lesssim \sum_{j=1}^{\infty}(2^j)^{\Gl+2n(1-s)-s\frac{\Gd}{\Gs}} (r^{\rho})^{\Gl+2n(1-s)+s\Gd(\frac{1}{\rho}-\frac{1}{\Gs})}r^{2ns(1-\frac{1}{p})} \\
    &\lesssim r^{\rho\Gl+2n[\rho(1-\frac{s}{2})+s(\frac{1}{q}-\frac{1}{p})]}\lesssim r^{\Gl+2n(1-\frac{s}{p})}
\end{align*}
provided that $\rho=\frac{2n(1-q^{-1})+\Gd}{n+\Gd\Gs^{-1}}$. Note that the infinite series converges by noting that 
$$\Gg_{p_0}=\frac{\Gg}{1-\rho}-n=\frac{\Gd}{\Gs}\frac{\Gs-\rho}{1-\rho}\leq \frac{\Gd}{\Gs},$$
and 
$$\frac{\Gl}{s}-\frac{2n}{s'}<\Gg_{p_0}\leq \frac{\Gd}{\Gs}.$$
% because  
% $-2n(1-\frac{s}{2})+\frac{2ns}{1-\rho}(\frac{1}{q}-\frac{1}{2})<2n(s-1)+s\Gd<n(s-1)+\frac{s\Gd}{\Gs}$ provided that 
% $$\Gl\leq -2n\bigg(1-\frac{s}{2}\bigg)+\frac{2ns}{1-\rho}\bigg(\frac{1}{q}-\frac{1}{2}\bigg)= -2n\bigg(1-\frac{s}{2}\bigg)+\frac{s\Gg(n+\frac{\Gd}{\Gs})}{\Gg+\frac{\Gd}{\Gs}-\Gd}.$$

Condition (M3') follows in the same way as in the proof of Theorem \ref{main1}. 

Thus, this completes the proof.
\end{proof}

% \subsection{A reduction on smoothness condition}
% \begin{lemma}
%     Under the assumption of \eqref{sizeK}, \eqref{smoothK} holds iff \eqref{smoothK} holds without $\Go(z,w)$. If one considers the non-convolution type with \eqref{sizeK'}, then 
%     \eqref{smoothK'} (or \eqref{smoothK'}) holds iff it holds without the factor $\Go(y,c)$
% \end{lemma}

\section{Dual spaces of $\sH^p$}\label{sec5}
Let us start with the definition of the twisted Campanato space.
\begin{definition}
    We say $f\in \GL^\times_{\Ga}(\C^n)$ if $f\in L^2_{loc}(\C^n)$ that %there exists $K>0$ such that for all ball $B=B(z_B,r)$, there exists a polynomial $P_{z_B}$ of degree at most $\lfloor\Ga\rfloor$ such that 
    $$\|f\|_{\GL^\times_{\Ga}}:=\sup_{r(B)<1}\inf_{P_{z_B}\in \sP_{\lfloor\Ga\rfloor}}\frac{1}{[r(B)]^{\Ga}}\bigg(\frac{1}{|B|}\int_B|f(w)-P_{z_B}(w)\Go(z_B,w)|^2dw\bigg)^{\frac12}+\sup_{r(B)\geq 1}\bigg(\frac{1}{|B|}\int_B|f(w)|^2dw\bigg)^{\frac{1}{2}},$$ 
    is finite, 
    where $\sP_{\ell}=\text{span}_{\C}\{z^{\Ga}\ol{z}^{\Gb}: |\Ga|+|\Gb|\leq \ell\}$.
    %The norm $\|f\|_{\GL^\times_{\Ga}}$ is defined to be the smallest possible $K$ for which the above holds.
\end{definition}
Note that $\GL^{\times}_0\simeq \mathscr{BMO}$ given in \cite{MPR}. For the sake of completeness, we give the proof of the duality here. %As in the classical case, we define $L^2_{N_p,0}(\C^n)$ be the space of $L^2_c(\C^n)$ that if $\supp(g)\subset B(c_B,R)$ with $R<1$, $$\int_{\C^n}(z-c_B)^{\Ga}\ol{(z-c_B)}^{\Gb}g(z)\Go(c_B,z)dz=0.$$

\begin{theorem} \label{dual1}
Let $0<p< 1$. Then,
\begin{enumerate}
    \item for every $g\in \GL^{\times}_{\Gg_p}$, the operator $L_g(f):=\int_{\C^n}f(z)g(z)dz$ initially defined on the $L^2_{c}(\C^n)$ can be extend uniquely to a bounded linear functional on $\sH^p$ and 
    $$\|L_g\|\leq c\|g\|_{\GL^{\times}_{\Gg_p}};$$
    \item conversely, every continuous linear functional $L$ on $\sH^p$ can be realized as a linear operator $L_g$ defined on $L^2_{c}(\C^n)$ with $g\in \GL^{\times}_{\Gg_p}$ and 
    $$\|g\|_{\GL^{\times}_{\Gg_p}}\leq C\|L\|.$$
\end{enumerate}
\end{theorem}
\begin{proof}\
\begin{enumerate}
    \item Suppose $a$ is a $(p,2)$-atom with $\supp(a)\subset B=B(z_B,r)$. Then, if $r<1$, 
    \begin{align*}
        \bigg|\int_{\C^n} a(z)g(z)dz\bigg|&=\bigg|\int_{\C^n} a(z)[g(z)-P_{\lfloor\Ga\rfloor,z_B}(z)\omega(z_B,z)]dz\bigg|\\
        &\leq \|a\|_{L^2}\|g-P_{N_p,z_B}\omega(z_B,\cdot)\|_{L^2(B)} \\
        &\leq \frac{1}{|B|^{\frac{1}{p}-\frac{1}{2}}}\|g-P_{N_p,z_B}\omega(z_B,\cdot)\|_{L^2(B)}.
    \end{align*}
    For $r\geq1$, simply add
\begin{align*}
    \bigg|\int ag\bigg|&\leq \|a\|_{L^2}\|g\|_{L^2(B)}\lesssim  r^{2n(1-1/p)}
\left(\fint_B|g|^2\right)^{1/2}\lesssim\|g\|_{\Lambda^\times_{\gamma_p}}.
\end{align*}
Here we used $r^{-\gamma_p}\leq1$.
    
    Now taking the infimum over all polynomials of degree at most $N_p$ on both sides, we have 
    \begin{align*}
       \bigg|\int_{\C^n} a(z)g(z)dz\bigg|&\leq \frac{c}{r^{n+\Gg_p}}\|g-P_{N_p,z_B}\omega(z_B,\cdot)\|_{L^2(B)} \\
        &\lesssim \sup_{B}\inf_{P_{z_B}}\frac{1}{r^{\Gg_p}}\|g-P_{N_p,z_B}\omega(z_B,\cdot)\|_{L^2(B;dx/|B|)}\simeq \|g\|_{\GL^{\times}_{\Gg_{p}}}.
    \end{align*}
    That is, $|L_g(a)|\lesssim \|g\|_{\GL^{\times}_{\Gg_{p}}}.$ Since the finite atomic decomposition is dense in $\sH^p$ with the aid of Theorem \ref{ft1}, we can extend $L_g$ uniquely to $\sH^p$ and $|L_g(f)|\leq c \|f\|_{\sH^p}\|g\|_{\GL^{\times}_{\Gg_{p}}}$ , which proves the first part.
   % Then, since $g$ is bounded, a standard argument shows that $|L_g(f)|\leq c \|f\|_{\sH^p}\|g\|_{\GL^{\times}_{\Gg_{p}}}$ for finite linear combination of $f$. Since finite atomic decomposition is dense in twisted Hardy spaces and by Theorem \ref{ft1}, we have a unique extension of $L$ to $\sH^p(\C^n)$ with $\|L\|\leq C'\|g\|_{\GL^{\times}_{\Gg_{p}}}$.

    \item Suppose $L:\sH^p(\C^n)\to \C$. Define $B_R:=B(0,R)\subset \C^n$. For any $f\in L^2(B_R)$ with $R\geq 1$, we can extend it to a function in $\sH^p(\C^n)$ by zero extension (and this would correspond to atoms with $r(B)>1$). Thus, $$|L(f)|\leq \|L\|\|f\|_{\sH^p}\leq C_n\|L\|R^{2n(\frac{1}{p}-\frac{1}{2})}\|f\|_{L^2(B_R)}.$$
    Thus, by Riesz representation theorem, there exists $g_R\in L^2(B_R)$ such that $L(f)=\int_{B_R} fg_R$ for $f\in L^2(B_R)$ and $\|g_R\|_{L^2(B_R)}\leq C\|L\| r^{2n(\frac{1}{p}-\frac{1}{2})}$. Since $B_R\subset B_{R+1}$, $f\in L^2(B_R)$ with zero extension is also in $L^2(B_{R+1})$. Thus, $g_{R}=g_{R+1}$ a.e. on $B_R$. Now, define $g(x)=g_R(x)$ if $x\in B_R$. Then $L(f)=\langle f,g\rangle$ for all $f\in L^2_c(\C^n)$. It remains to prove $g\in \GL^\times_{\Ga}$.

    If $r(B)\geq1$, then $\sup_{c\in \C^n}\|g\|_{L^2(B(c,1))}\leq C\|L\|$ by taking $a=\frac{f}{\|f\|_{L^2(B(c,1))}}$ and then taking supremum over all $f\in L^2(B(c,1))$. For any balls $B(c,r)$ with $r>1$, we can cover $B(c,r)$ by $O(r^{2n})$ many balls of radius 1. Then, 
    $$\bigg(\fint_{B}|g(z)|^2dz\bigg)^{\frac{1}{2}}\lesssim \frac{r^{n}}{|B|^{\frac{1}{2}}}C\|L\|\lesssim \|L\|.$$
    %and this implies $$\sup_{r(B)\geq  1}\bigg(\frac{1}{[r(B)]^{2n}}\int_{B}|g(z)|^2dz\bigg)^{\frac{1}{2}}\leq C\|L\|.$$

    If $r(B)<1$, consider the closed subspace 
  $$
V_B=
\left\{
f\in L^2(B):
\int_B f(z)P(z)\Omega(c_B,z)\,dz=0,\
\forall P\in\mathscr P_{N_p}
\right\}.
$$
Now since elements in $V_B$ are also elements in $\sH^p(\C^n)$, we have 
    $$\inf_{P\in \sP_{N_p}}\|g-P\Go(c_B,\cdot)\|_{L^2(B)}=\sup_{f\in V_B, \|f\|_{L^2}\leq 1}\bigg|\int fg\bigg|\leq C\|L\|r^{2n(\frac{1}{p}-\frac{1}{2})},$$
    where $P$ is a polynomial of degree up to $N_p$. Thus,
    $$\frac{1}{r^{2n(1/p-1)}}\bigg(\frac{1}{r^{n}}\inf_{P\in \mathscr{P}_{N_p}}\|g-P\Go(c_B,\cdot)\|_{L^2(B)}\bigg)\leq C\|L\|.$$
    Thus, $g\in \GL^{\times}_{\Gg_p}$.
\end{enumerate}
\end{proof}

Next, we would like to characterize twisted Campanato spaces using the twisted derivative. 
\begin{definition} \label{GLxat}
    Let $\Ga>0$. If $\Ga\notin \N$, we define 
    $$\widetilde{\GL}^{\times}_{\Ga}:=\bigg\{g\in C^{\lfloor\Ga\rfloor}(\C^n): [g]_{\widetilde{\GL}^{\times}_{\Ga},J}:=\sup_{z\in \C^n, |h|\in (0,1)} \frac{|\tau^{\widetilde{X}}_{h}(\widetilde{X}^Jg)(z)-\widetilde{X}^Jg(z)|}{|h|^{\Ga-\lfloor\Ga\rfloor}}<\infty \quad \forall\ J\in \{1,\cdots,2n\}^{\lfloor\alpha\rfloor}\bigg\};$$
    and if $\Ga\in \N$, we define 
    \begin{align*}
        \widetilde{\GL}^{\times}_{\Ga}:=\bigg\{g\in C^{\Ga-1}(\C^n): [g]_{\widetilde{\GL}^{\times}_{\Ga},J}:=&\sup_{z\in \C^n, |h|\in (0,1)} \frac{|\tau^{\widetilde{X}}_{h}(\widetilde{X}^Jg)(z)-2\widetilde{X}^Jg(z)+\tau^{\widetilde{X}}_{-h}(\widetilde{X}^Jg)(z)|}{|h|}<\infty,\\
        &\quad \quad \quad \quad \forall J\in \{1,\cdots,2n\}^{\Ga-1}\bigg\}.
    \end{align*}

    The norm of $\widetilde{\GL}^{\times}_{\Ga}$ is given by 
$$\|g\|_{\widetilde{\GL}^{\times }_{\Ga}}:=\sum_{J}[g]_{\widetilde{\GL}^{\times}_{\Ga},J}+\|g\|_{L^{\infty}},$$
where the sum is taken over all multiindices $J=(j_1,\cdots, j_{|J|})\in \{1,\cdots, 2n\}^{|J|}$.
\end{definition}

% The norm of $\widetilde{\GL}^{\times}_{\Ga}$ is given by 
% $$\|g\|_{\widetilde{\GL}^{\times }_{\Ga}}:=\sum_{J}[g]_{\widetilde{\GL}^{\times}_{\Ga},J}+\|g\|_{L^{\infty}},$$
% where the sum is taken over all multiindices $J=(j_1,\cdots, j_{|J|})\in \{1,\cdots, 2n\}^{|J|}$.

Pick a $\eta\in C^{\infty}_c(\C^n)$ with $\supp(\eta)\subset B(0,1)$, $\int\eta =1$ and $\int \eta(z)z^{\Ga}\ol{z}^{\Gb}dz=0$ with order up to $1\leq |\Ga|+|\Gb|\leq N_p+k$. Consider 
$a_{c,t}(z)=\ds Ct^{2n+k+\ell-\frac{2n}{p}} \frac{d^\ell}{dt^\ell}[\widetilde{X}^{(i_1,\cdots,i_k)}\eta_t(c-z)]\Go(z,c)$.

\begin{lemma} \label{act-atom}
For $k\in \N\cup\{0\}$,  $\ell\in \N$, $c\in \C^n$, and $t\in (0,1)$, the function $a_{c,t}$ is a constant multiple (independent of $t$) of an $(p,\infty)$-atom, where the constant is independent of $c$ and $t$, provided that $\eta$ has vanishing moment conditions up to order $N_p+k$.
\end{lemma}
\begin{proof}
    
It can be seen that $\supp(a_{c,t})\subset B(c,t)$. Moreover, since
    $$\|\partial_t^\ell\widetilde{X}^{(i_1,\cdots,i_k)}\eta_t\|_{L^\infty}\lesssim t^{-2n-k-\ell},$$
    we have$\|a_{c,t}\|_{L^\infty}\leq C_{n} t^{-\frac{2n}{p}}$ whenever $t\in (0,1)$. 
    
     For the cancellation condition, let $|\Ga|+|\Gb|\leq N_p$. By a change of variables and integration by parts, we have
\begin{align*}
 &\int_{\C^n} a_{c,t}(w)(w-c)^{\Ga}\ol{(w-c)}^{\Gb}\Go(c,w)dw\\
 &=Ct^{2n+k+\ell-\frac{2n}{p}}\frac{d^\ell}{dt^\ell}\int_{\C^n}  \widetilde{X}^{(i_1,\cdots,i_k)}\eta_t(c-w) (w-c)^{\Ga}\ol{(w-c)}^{\Gb}dw\\
 &= Ct^{2n+k+\ell-\frac{2n}{p}} (-1)^{|\Ga|+|\Gb|}\frac{d^\ell}{dt^\ell}\int_{\C^n}   \widetilde{X}^{(i_1,\cdots,i_k)} \eta_t(w)w^{\Ga}\ol{w}^{\Gb}dw\\
% &= (-1)^{|\Ga|+|\Gb|}\frac{d^\ell}{dt^\ell}\int_{\C^n}  \widetilde{X}^{(i_1,\cdots,i_k)}\ol{\eta_t}(w)w^{\Ga}\ol{w}^{\Gb}dw\\
 &= Ct^{2n+k+\ell-\frac{2n}{p}} (-1)^{|\Ga|+|\Gb|+k}\frac{d^\ell}{dt^\ell}\int_{\C^n} \eta_t(w)X^{(i_k,\cdots,i_1)}[w^{\Ga}\ol{w}^{\Gb}]dw
\end{align*}
Now we know that from the definition of $\widetilde{X}$, we have 
$$X^{(i_k,\cdots,i_1)}(w^{\Ga}\ol{w}^{\Gb})=\sum_{|\Gg|+|\delta|\leq k+N_p} b_{\Gg,\Gd} w^{\Gg}\ol{w}^{\Gd}.$$
Thus, by the cancellation conditions of $\eta$, we obtain
$$\int_{\C^n}\eta_t(w)X^{(i_k,\cdots,i_1)}  (w^\Ga\ol{w}^\Gb)dw=b_{0,0}.$$
Since $\ell\geq 1$, the derivative with respect to $t$ yields the term 0. 

Therefore, $a_{c,t}$ is a constant multiple of a $(p,\infty)$-atom.
\end{proof}

\begin{theorem} \label{e1}
    Suppose $\alpha>0$. Then, $\GL^{\times}_{\Ga}\hookrightarrow \widetilde{\GL}^{\times}_{\Ga}$.
\end{theorem}

\begin{proof}
We follow the idea of \cite{WeightedIneq}*{Chapter III Theorem 5.22}
    Let $g\in {\GL}^{\times}_{\Ga}$. Pick $\phi\in C^{\infty}_c(\R^{2n})$ with $\supp(\phi)\subset B(0,1)$, $\int \phi=1$ and $\int_{\C^n}\phi(z)z^{\Gd}\ol{z}^{\Gb}=0$ for $1\leq |\Gd|+|\Gb|\leq 2N_p+2$. Note that $\{\phi_t\}_{t>0}$ forms an approximation to the identity under both standard convolution and twisted convolution.

    Define $u_g(z,t)=\phi_t\times g(z)$. Then pick $p\in (0,1]$ such that $\Ga=2n(p^{-1}-1)$ and $k=N_p$. Then,
    $$|\partial^\ell_t\widetilde{X}_{(i_1,\cdots, i_m)}u_g(z,t)|=|(\partial^\ell_t\widetilde{X}_{(i_1,\cdots, i_m)}\phi_t)\times g(z)|\leq Ct^{\Ga-m-\ell}\|g\|_{\GL^{\times}_{\Ga}},$$
    where $C$ is independent of $z$, $t$ and $g$. Here, $\ell\geq1$ and $0\leq m\leq N_p+2$.

    We first show that $g\in L^{\infty}(\C^n)$. Taking $m=0$ and $\ell=1$,
    $|\partial_tu_g(z,t)|\lesssim t^{\Ga-1}\|g\|_{\GL^\times_{\Ga}},$
    and $|u_g(z,1)|\leq \|\phi\|_{L^2}\|g\|_{\GL^\times _{\Ga}}\lesssim\|g\|_{\GL^\times _{\Ga}},$
    thus, 
    \begin{align*}
        |u_g(z,t)|&\leq |u_g(z,1)|+\int_t^1 |\partial_tu_g(z,t)|dt \lesssim \|g\|_{\GL^\times_{\Ga}},
    \end{align*}
    and the implicit constant is independent of $t$. Letting $t\to0^+$, one has 
    $$|g(z)|\lesssim \|g\|_{\GL^{\times }_{\Ga}}.$$

    Next, we consider $\Ga\notin \N$. Then, $\Ga=N_p+\Ga'$, where $\Ga'\in (0,1)$. Given $|z-c_B|<1$, if $\ell=1$, we have 
    \begin{align*}
        |\widetilde{X}_{(i_1,\cdots, i_k)}u(z,t_2)-\widetilde{X}_{(i_1,\cdots, i_k)}u(z,t_1)|&\leq C\|g\|_{\GL^{\times}_{\Ga}} \int_{t_1}^{t_2}t^{\Ga'-1}dt\leq C\|g\|_{\GL^{\times}_{\Ga}} t_2^{\Ga'}.
    \end{align*}
    Thus, as $t\to 0^+$, the function $z\mapsto \widetilde{X}_{(i_1,\cdots, i_k)}u(z,t)$ will converge uniformly to a continuous function $u_{i_1,\cdots, i_k}(z)$. Meanwhile, the approximation to identity gives us 
    \begin{align*}
        \langle \widetilde{X}_{(i_1,\cdots,i_k)}(\phi_t\times g), \psi\rangle &= \langle (\phi_t\times g), [\widetilde{X}_{(i_1,\cdots,i_k)}]^*\psi\rangle\\
        &\to \langle g, [\widetilde{X}_{(i_1,\cdots,i_k)}]^*\psi\rangle=\langle\widetilde{X}_{(i_1,\cdots,i_k)}g,\psi\rangle
    \end{align*}
    as $t\to 0^+$ for all $\psi\in C^{\infty}_c(\C^n)$.
    Thus, uniqueness of limit gives us $\widetilde{X}_{(i_1,\cdots,i_k)}g$ is continuous. The same argument applies to all lower-order twisted derivatives. Hence, by Remark \ref{rem:regularity}, we have $g\in C^k(\C^n)$.

    Now for $z\in \C^n$ and $|h|<1$, we have 
    \begin{align*}
        &|\tau^{\widetilde{X}}_{h}u_{i_1,\cdots, i_k}(z)-u_{i_1,\cdots, i_k}(z)|\\&\leq |\tau^{\widetilde{X}}_{h}[u_{i_1,\cdots, i_k}(z)-\widetilde{X}_{(i_1,\cdots,i_k)}(\phi_{|h|}\times g)(z)]|+|u_{i_1,\cdots, i_k}(z)-\widetilde{X}_{(i_1,\cdots,i_k)}(\phi_{|h|}\times g)(z)|\\
        &\quad\quad +|\tau^{\widetilde{X}}_{h}\widetilde{X}_{(i_1,\cdots,i_k)}(\phi_{|h|}\times g)(z)-\widetilde{X}_{(i_1,\cdots,i_k)}(\phi_{|h|}\times g)(z)|
    \end{align*}
    The first two terms give the upper bound $C\|g\|_{\GL^{\times}_{\Ga}} |h|^{\Ga'}$. The last term we need to deal with care. First note that 
    \begin{align*}
        \widetilde{X}_{(i_1,\cdots,i_k)}(\phi_{|h|}\times g)(z) &=-\int_{|h|}^1\partial_s\widetilde{X}_{(i_1,\cdots,i_k)}(\phi_{s}\times g)(z)ds+\widetilde{X}_{(i_1,\cdots,i_k)}(\phi_{1}\times g)(z).
    \end{align*}
The term involving $\phi_1\times g$ is easier because (from Taylor expansion)
\begin{align*}
    |(\tau^{\widetilde{X}}_{h}-1)\widetilde{X}_{(i_1,\cdots,i_k)}(\phi_{1}\times g)(z)|&\leq |h| \sum_{j=1}^{2n}\sup_{w\in B(z,|h|)}|(\widetilde{X}_{(j,i_1,\cdots,i_k)}\phi_1)\times g(w)| \\
    &\lesssim|h| \|g\|_{L^{\infty}}.
\end{align*}
Thus, as $|h|<1$, we have 
$$|(\tau^{\widetilde{X}}_{h}-1)\widetilde{X}_{(i_1,\cdots,i_k)}(\phi_{1}\times g)(z)|\leq C|h|^{\Ga'} \|g\|_{L^{\infty}}.$$

At the same time, we have 
\begin{align*}
&\bigg|\int_{|h|}^1(\tau_h^{\widetilde{X}}-1)\partial_s\widetilde{X}_{(i_1,\cdots,i_k)}(\phi_{s}\times g)(z)ds\bigg|\\
    &\leq \int_{|h|}^1 |h|\sum_{j=1}^{2n}\sup_{w\in B(z,|h|)}|(\partial_s\widetilde{X}_{j,i_1,\cdots,i_k}\phi_s)\times g(w)| ds \\
    &\lesssim |h|\int_{|h|}^1 s^{\Ga-N_p-2}\|g\|_{\GL^{\times }_{\Ga}} ds \lesssim |h||h|^{\Ga-N_p-1}\|g\|_{\GL^{\times }_{\Ga}} =|h|^{\Ga'}\|g\|_{\GL^{\times }_{\Ga}} .
\end{align*}
Thus, we can conclude that 
$$\sup_{z\in \C^n, |h|<1} \frac{|\tau^{\widetilde{X}}_{h}\widetilde{X}_{(i_1,\cdots,i_k)}g(z)-\widetilde{X}_{(i_1,\cdots,i_k)}g(z)|}{|h|^{\Ga'}}\leq C'\|g\|_{\GL^{\times }_{\Ga}}$$
for any $(i_1,\cdots,i_k)\in \{1,\cdots,2n\}^{k}$.

When $\Ga\in \N$, Consider $k=N_p-1$, then $\Ga<k+2$ and for $t\in (0,1)$
$$|\partial_t^2\widetilde{X}_{(i_1,\cdots,i_{k})}u_g(z,t)|\leq C\|g\|_{\GL^{\times }_{\Ga}}t^{-1}.$$
Then, 
$$|\partial_t\widetilde{X}_{(i_1,\cdots,i_{k})}u_g(z,t)|\leq\Big|\partial_t\widetilde{X}_{(i_1,\cdots,i_{k})}u_g(z,t)|_{t=1}\Big|+C\|g\|_{\GL^{\times }_{\Ga}}\ln(t^{-1}).$$
Now note that 
$$\partial_t\widetilde{X}_{(i_1,\cdots,i_{k})}u_g(z,t)|_{t=1}=(\partial_t\widetilde{X}_{(i_1,\cdots,i_{k})}\phi_t|_{t=1})\times g(z),$$
by Lemma \ref{act-atom}, we have 
$$\Big|\partial_t\widetilde{X}_{(i_1,\cdots,i_{k})}u_g(z,t)|_{t=1}\Big|\leq C_{\phi}\|g\|_{L^\infty}.$$
Thus, 
$$|\partial_t\widetilde{X}_{(i_1,\cdots,i_{k})}u_g(z,t)|\leq C\|g\|_{\GL^{\times }_{\Ga}}\ln(et^{-1}).$$
Here, $C$ is independent of $z$, $t$ and $g$. For $0<t_1<t_2<1$, we have 
$$|\widetilde{X}_{(i_1,\cdots,i_{k})}u_g(z,t_2)-\widetilde{X}_{(i_1,\cdots,i_{k})}u_g(z,t_1)|\leq C'\|g\|_{\GL^{\times }_{\Ga}}t_2\ln(et_2^{-1}),$$
which approaches 0 as $t_2\to 0^+$. Therefore, following the same argument as in the case $\Ga\notin \N$, the limit is $\widetilde{X}_{(i_1,\cdots,i_{k})}g(z)$ and is continuous. Again, the same argument applies to all lower-order twisted derivatives. Hence, by Remark \ref{rem:regularity}, we have $g\in C^k(\C^n)$.

Moreover, since $t\ln(et^{-1})\to 0$ as $t\to 0^+$,
we have $$t\partial_t\widetilde{X}_{(i_1,\cdots,i_k)}u_g(z,t)\to 0$$
uniformly in $z$. Hence, by integration by parts,
we can write
\begin{align*}
    &\widetilde{X}_{(i_1,\cdots,i_{k})}g(z)\\
    &=\widetilde{X}_{(i_1,\cdots,i_{k})}u_g(z,|h|)-|h|\partial_t\widetilde{X}_{(i_1,\cdots,i_{k})}u_g(z,|h|)+\int_0^{|h|} s \partial^2_s\widetilde{X}_{(i_1,\cdots,i_{k})}u_g(z,s)ds=I+II+III.
\end{align*}
Then, a direct estimate gives us
\begin{align*}
    &|\tau^{\widetilde{X}}_hIII-2III+\tau^{\widetilde{X}}_{-h}III|\leq 4C\int_0^{|h|}s s^{-1}ds \|g\|_{\GL^{\times }_{\Ga}}=C'\|g\|_{\GL^{\times }_{\Ga}} |h|
\end{align*}
and 
\begin{align*}
    &|\tau^{\widetilde{X}}_hII-2II+\tau^{\widetilde{X}}_{-h}II|\lesssim|h|^3 \sum_{i,j=1}^{2n}\sup_{w\in B(z,|h|)}|\widetilde{X}_{(i,j)}\partial_t\widetilde{X}_{(i_1,\cdots,i_{k})}u_g(w,|h|)|\lesssim_{\phi} \|g\|_{\GL^{\times }_{\Ga}}|h|.
\end{align*}

It remains to estimate $I$. Note that 
\begin{align*}
    \widetilde{X}_{(i_1,\cdots,i_{k})}u_g(z,|h|)&= \widetilde{X}_{(i_1,\cdots,i_{k})}u_g(z,1)-\int_{|h|}^{1} \partial_s \widetilde{X}_{(i_1,\cdots,i_{k})}u_g(z,s)ds.
\end{align*}
Thus,
\begin{align*}
    &(\tau^{\widetilde{X}}_h-2I+\tau^{\widetilde{X}}_{-h})\widetilde{X}_{(i_1,\cdots,i_{k})}u_g(z,|h|) \\
    &=(\tau^{\widetilde{X}}_h-2I+\tau^{\widetilde{X}}_{-h})\widetilde{X}_{(i_1,\cdots,i_{k})}u_g(z,1)-\int_{|h|}^{1} (\tau^{\widetilde{X}}_h-2I+\tau^{\widetilde{X}}_{-h}) \partial_s \widetilde{X}_{(i_1,\cdots,i_{k})}u_g(z,s)ds
\end{align*}
and thus
$$|(\tau^{\widetilde{X}}_h-2I+\tau^{\widetilde{X}}_{-h})\widetilde{X}_{(i_1,\cdots,i_{k})}u_g(z,1)|\lesssim |h|^2\|g\|_{\GL^{\times }_{\Ga}}\lesssim |h|\|g\|_{\GL^{\times }_{\Ga}}$$
and 
\begin{align*}
    &\bigg|\int_{|h|}^{1} (\tau^{\widetilde{X}}_h-2I+\tau^{\widetilde{X}}_{-h}) \partial_s \widetilde{X}_{(i_1,\cdots,i_{k})}u_g(z,s)ds\bigg|\lesssim|h|^2\int_{|h|}^1 \|g\|_{\GL^{\times}_{\Ga}}s^{-2}ds \simeq |h|\|g\|_{\GL^{\times}_{\Ga}}.
\end{align*}

Combining all these estimates, we have 
$$\sup_{z\in \C^n, |h|<1} \frac{|\tau^{\widetilde{X}}_{h}\widetilde{X}_{(i_1,\cdots,i_k)}g(z)-2\widetilde{X}_{(i_1,\cdots,i_k)}g(z)+\tau^{\widetilde{X}}_{-h}\widetilde{X}_{(i_1,\cdots,i_k)}g(z)|}{|h|}\leq C'\|g\|_{\GL^{\times }_{\Ga}}$$
for any $(i_1,\cdots,i_k)\in \{1,\cdots,2n\}^{k}$.
\end{proof}

To show the converse, we first consider the easier case, $\alpha \notin \N$.

\begin{lemma}
    Let $\Ga>0$ and $\Ga\notin \N$. If $g\in \widetilde{\GL}^{\times}_{\Ga}$, then $g\in \GL^{\times}_{\Ga}$.
\end{lemma}
\begin{proof}
We first note that the case $N=0$ is immediate, because 
$|g(z)-g(c_B)\Go(c_B,z)|\leq |z-c_B|^{\Ga}\|g\|_{\widetilde{\GL}^\times_{\Ga}}$ and $g$ is bounded.

   The Taylor polynomial of degree $N:=\lfloor \Ga\rfloor>0$ at $c_B$ of $g$ is given by  (applying \eqref{rightT} with $z-z_B=-c_B$, $w-z_B=z-c_B$ and replacing $g(-\cdot)$ by $g$ and absorb $(-1)^k$)
    \begin{align*}
    g(z)&=g(c_B)\Go(c_B,z)+\sum_{k=1}^{N-1}\sum_{1\leq i_1,\cdots, i_k\leq 2n} \frac{1}{k!} T_{i_1,\cdots, i_k}(z-c_B)[g(c_B)]\Go(c_B,z)\\
    &\quad +\sum_{i_1,\cdots,i_{N}=1}^{2n}\int_0^1\frac{(1-\Gth)^{N-1}}{(N-1)!}e^{-\frac{i}{2}\Gth\Im(c_B\cdot\ol{z})}T_{i_1,\cdots,i_{N}}(z-c_B)[g(\Gth (z-c_B)+c_B)]d\Gth \Go(c_B,z)
\end{align*}
If we define 
$$P(z)=g(c_B)+\sum_{k=1}^{N}\sum_{1\leq i_1,\cdots, i_k\leq 2n} \frac{1}{k!} T_{i_1,\cdots, i_k}(z-c_B)[g(c_B)]$$
and expand $h$ up to $(N-1)$-th degree with error in degree $N$, we have 
\begin{align*}
    &g(z)-P(z)\Go(c_B,z)\\
    &= \Go(c_B,z)\sum_{i_1,\cdots,i_{N}=1}^{2n}\int_0^1\frac{(1-\Gth)^{N-1}}{(N-1)!}e^{\frac{i}{2}\theta\Im((-c_B)\cdot\ol{z-c_B})}T_{i_1,\cdots,i_{N}}(z-c_B)[g(\theta(z-c_B)+c_B)]d\Gth\\
    &\quad\quad - \Go(c_B,z)\sum_{i_1,\cdots,i_{N}=1}^{2n}\int_0^1\frac{(1-\Gth)^{N-1}}{(N-1)!}T_{i_1,\cdots,i_{N}}(z-c_B)[g(c_B)]d\Gth.
\end{align*}
Then, since $|z-c_B|<r<1$, by writing $\Ga=N+\Ga'$,
\begin{align*}
    &|g(z)-P(z)\Go(c_B,z)|\\
    &\leq |z-c_B|^{N} \sum_{i_1,\cdots,i_{N}=1}^{2n}\int_0^1\frac{(1-\theta)^{N-1}}{(N-1)!}|\tau^{\widetilde{X}}_{-\theta(z-c_B)}([\widetilde{X}_{(i_1,\cdots, i_N)}g](c_B))-\widetilde{X}_{(i_1,\cdots, i_N)}g(c_B)|d\theta \\
    &\lesssim |z-c_B|^{N} \|g\|_{\widetilde{\GL}^{\times}_{\Ga}}\sum_{i_1,\cdots,i_{N},j=1}^{2n}\int_0^1\frac{(1-\theta)^{N-1}}{(N-1)!}|\theta(z-c_B)|^{\Ga'}d\theta \\
    &\lesssim |z-c_B|^{\Ga}\|g\|_{\widetilde{\GL}^{\times}_{\Ga}}.
\end{align*}
Thus, if $0<r<1$
$$\bigg(\frac{1}{|B(c_B,r)|}\int_{B(c_B,r)}|g(z)-P(z)\Go(c_B,z)|^2dz\bigg)^{\frac12}\lesssim r^{\Ga} \|g\|_{\widetilde{\GL}^{\times}_{\Ga}}, $$
and the boundedness of $g$ implies $\frac{1}{r^{2n}}\int_{B(x,r)}|g(z)|^2dz<\infty$ for all $r\geq 1$, which shows $g\in \GL^{\times }_{\Ga}$.
\end{proof}

For $\Ga\in \N$, we cannot follow the proof in  \cite{WeightedIneq}, because the twisted derivative will not necessarily lower the degree of a polynomial. The idea is to make use of ``finite difference Zygmund space". More precisely, we consider the space $g\in L^\infty(\C^n)$ with 
$$\|g\|_{\widetilde{\GL'}^{\times}_{\Ga}}:=\|g\|_{L^\infty}+\sup_{|h|<1}|h|^{-\alpha}\|(\tau^{\widetilde{X}}_h-I)^{\alpha+1}g\|_{L^\infty}<\infty.$$

Then we aim to transfer the problem to the classical Zygmund and Campanato spaces, and then transfer the result back to the twisted case.

\begin{lemma} \label{e2}
    Let $\Ga\in \N$. For $f\in \widetilde{\GL}^{\times}_{\Ga}$, we have 
   % $$\sup_{|h|<1, z\in\C^n}\frac{|(\tau_h^{X}-I)^{\lceil\Ga\rceil}f(z)|}{|h|^{\Ga}}\leq C\|f\|_{\widetilde{\GL}^{\times}_{\Ga}}$$
    %if $\Ga\notin \N$ and 
    $$\sup_{|h|<1, z\in\C^n}\frac{|(\tau_h^{\widetilde{X}}-I)^{\Ga+1}f(z)|}{|h|^{\Ga}}\leq C\|f\|_{\widetilde{\GL}^{\times}_{\Ga}},$$
    where the constant $C$ is independent of $f$.

    That is, $\|f\|_{\widetilde{\GL'}^{\times}_{\Ga}}\leq C\|f\|_{\widetilde{\GL}^{\times}_{\Ga}}$.
\end{lemma}
\begin{proof}
This is a twisted version of the classical result. Define the directional derivative 
    $$D_hf:=\frac{d}{dt}(\tau^{\widetilde{X}}_{th}f)\bigg|_{t=0}=-h\cdot \nabla f(z)+\frac{i}{2}\Im(z\cdot \ol{h})f(z).$$
    Its operator form is given by 
    $$D_h=-\sum_{j=1}^{2n}h_j \widetilde{X}_j$$
    and its higher-order is 
    $$D_h^k=(-1)^k\sum_{j_1,\cdots,j_k=1}^{2n}h_{j_1}\cdots h_{j_k}\widetilde{X}_{(j_1,\cdots,j_k)}.$$
    Then,
    $$(\tau^{\widetilde{X}}_{h}-1)f=\int_0^1\frac{d}{dt}\tau^{\widetilde{X}}_{th}fdt$$
    and 
    \begin{align*}
        \frac{d}{dt}[f(z-th)e^{\frac{i}{2}\Im(z\cdot t\ol{h})}]&=-h\cdot\nabla f(z-th)e^{\frac{i}{2}\Im(z\cdot t\ol{h})}+\frac{i}{2}\Im(z\cdot \ol{h})e^{\frac{i}{2}\Im(z\cdot t\ol{h})}f(z-th)=\tau^{\widetilde{X}}_{th}D_hf(z).
    \end{align*}
    Indeed, we have 
    $$(\tau^{\widetilde{X}}_h-1)^{k}f(z)=\int_{[0,1]^k}\tau^{\widetilde{X}}_{(t_1+\cdots+t_k)h}D_h^kf(z)dt_1\cdots dt_k.$$
    Here we have used the fact that 
$$\tau^{\widetilde{X}}_{sh}\tau^{\widetilde{X}}_{th}=\tau^{\widetilde{X}}_{(s+t)h}=\tau^{\widetilde{X}}_{th}\tau^{\widetilde{X}}_{sh}$$
for $s,t\in \R$ and $h\in \C^n$.

    % If $\Ga\in (k-1,k)$, then 
    % \begin{align*}
    %     |(\tau^X_h-1)^kf(z)|&=|(\tau^X_h-1)[(\tau^X_h-1)^{k-1}f](z)|\\
    %     &\leq \bigg|\sum_{j_1,\cdots,j_{k-1}=1}^{2n}h_{j_1}\cdots h_{j_k}\int_{[0,1]^{k=1}}\tau^X_{(t_1+\cdots+t_{k-1})h}(\tau^X_h-1)\widetilde{X}_{(j_1,\cdots,j_{k-1})}f(z)dt_1\cdots dt_{k-1}\bigg|\\
    %     &\leq |h|^{k-1}\|(\tau^X_h-1)\widetilde{X}_{(j_1,\cdots,j_{k-1})}f\|_{L^{\infty}(\C^n)}\lesssim |h|^{\Ga}\|f\|_{\widetilde{\GL}^{\times}_{\Ga}}.
    % \end{align*}

    If $\Ga=k$, 
    \begin{align*}
        |(\tau^{\widetilde{X}}_h-1)^{k+1}f(z)|&=|(\tau^{\widetilde{X}}_h-1)^{k-1}[(\tau^{\widetilde{X}}_h-1)^{2}f](z)|\\
        &\leq |h|^{k-1}\sum_{j_1,\cdots,j_{k-1}=1}^{2n}\|(\tau^{\widetilde{X}}_{2h}-2\tau^{\widetilde{X}}_h+I)\widetilde{X}_{(j_1,\cdots, j_{k-1})}f\|_{L^\infty(\C^n)}\\
        &\lesssim_{\alpha}|h|^{k-1}\sum_{j_1,\cdots,j_{k-1}=1}^{2n}\|(\tau^{\widetilde{X}}_h-2I+\tau^{\widetilde{X}}_{-h})\widetilde{X}_{(j_1,\cdots, j_{k-1})}f\|_{L^\infty(\C^n)}\simeq |h|^{k}\|f\|_{\widetilde{\GL}^{\times}_{\Ga}}.
    \end{align*}
\end{proof}

Define $F_c(z)=\tau^{X}_cg(z)$.
\begin{lemma}\label{lemma1}
    Let $\Delta_hf(z):=f(z+h)-f(z)$, $m\in \N$, and $R>1$. Then 
    $$\Delta^m_hF_c(z-c)=\Go(z,c)(\Go(u,h)\tau^{\widetilde{X}}_{-h}-1)^{m}g(z)$$
    and 
    $$|\Delta^m_hF_c(z-c)|\leq C_{m,R}\sum_{\ell=0}^m|h|^{m-\ell}\|(\tau^{\widetilde{X}}_{-h}-I)^{\ell}g\|_{L^{\infty}}$$
    for all $z-c, z-c+mh\in B(0,R)$.
\end{lemma}

\begin{proof}
    For simplicity, write $u=z-c$. For any integer $j\geq 0$, $F_c(u+jh)=g(z+jh)\Go(u+jh,c)$ while 
    $\tau^{\widetilde{X}}_{-jh}g(z)=g(z+jh)\Go(z,-jh)$. Thus, 
    $$F_c(u+jh)=\Go(u+jh,c)\Go(z,jh)\tau^{\widetilde{X}}_{-jh}g(z).$$

    Thus, 
    \begin{align*}
        \Delta^m_hF_c(u)&=\sum_{j=0}^m (-1)^{m-j} {m\choose j} F_c(u+jh) \\
        &= \sum_{j=0}^m (-1)^{m-j} {m\choose j} \Go(u+jh,c)\Go(z,jh)[\tau^{\widetilde{X}}_{-h}]^jg(z)\\
        &=\Go(u,c)\sum_{j=0}^m (-1)^{m-j} {m\choose j} [\Go(u,h)]^j[\tau^{\widetilde{X}}_{-h}]^jg(z)\\
        &=\Go(z,c) [\Go(u,h)\tau^{\widetilde{X}}_{-h}-I]^mg(z).
\end{align*}

To prove the second inequality, by writing $\Go(u,h)\tau^{\widetilde{X}}_{-h}-I=[\Go(u,h)-1]+\Go(u,h)[\tau^{\widetilde{X}}_{-h}-I]$, we have 
\begin{align*}
|\Delta_h^mF_c(u)|&\leq \sum_{\ell=0}^m {m\choose \ell}
|\Go(u,h)-1|^{m-\ell}
\|(\tau^{\widetilde{X}}_{-h}-I)^\ell g\|_\infty\\
&\leq  \sum_{\ell=0}^m {m\choose \ell}
(\frac{|u|}{2}|h|)^{m-\ell}
\|(\tau^{\widetilde{X}}_{-h}-I)^\ell g\|_\infty\\
&\lesssim_{R,m}  \sum_{\ell=0}^m {m\choose \ell}|h|^{m-\ell}
\|(\tau^{\widetilde{X}}_{-h}-I)^\ell g\|_\infty,
\end{align*}
where the implicit constant is $(R/2)^{m}$.
\end{proof}

The next step is to control the term $\|(\tau^{\widetilde{X}}_{-h}-I)^{\ell}g\|_{L^\infty}$ merely using $\|g\|_{\widetilde{\GL'}^{\times}_{\Ga}}$, where $\Ga\geq \ell$.

\begin{lemma}[Marchaud's inequality]
    Suppose $g\in \widetilde{\GL'}^{\times}_{\Ga}$, $\Ga\in \N$, and $|h|<1$, then 
    $$\|(\tau^{\widetilde{X}}_{-h}-I)^\ell g\|_\infty
\leq
\begin{cases}
C_\Ga \|g\|_{{\widetilde{\GL'}^{\times}_{\Ga}}}|h|^\ell,&0\leq\ell<\Ga,\\
C'_\Ga \|g\|_{{\widetilde{\GL'}^{\times}_{\Ga}}}|h|^\Ga\left(1+\log\frac1{|h|}\right),
&\ell=\Ga,\\
\|g\|_{{\widetilde{\GL'}^{\times}_{\Ga}}} |h|^\Ga,&\ell=\Ga+1.
\end{cases}$$
\end{lemma}
Marchaud's inequality in the ordinary case can be found in many textbooks, for instance, \cite{DL93}*{Chapter 2 Section 8}.

\begin{proof}
    It is trivial for the case $\ell=\Ga+1$ or $|h|\in [\frac{1}{2},1)$. In the remainder of the proof, we will only consider $\ell\leq \Ga$ and $|h|<\frac{1}{2}$. 

    First note that 
    $$(\tau^{\widetilde{X}}_{-2h}-I)=(\tau^{\widetilde{X}}_{-h}-I)^2+2(\tau^{\widetilde{X}}_{-h}-I)$$
    and thus, for $1\leq \ell\leq \Ga$, we have 
   \begin{align*}
       (\tau^{\widetilde{X}}_{-2h}-I)^{\ell}&=\sum_{j=0}^{\ell} {\ell\choose j} 2^{\ell-j}(\tau^{\widetilde{X}}_{-h}-I)^{\ell+j}\\ 
       &= 2^{\ell}(\tau^{\widetilde{X}}_{-h}-I)^{\ell}+\sum_{j=1}^{\ell} {\ell\choose j} 2^{\ell-j}(\tau^{\widetilde{X}}_{-h}-I)^{\ell+j},
   \end{align*}
   which implies
   $$(\tau^{\widetilde{X}}_{-h}-I)^{\ell}=2^{-\ell} (\tau^{\widetilde{X}}_{-2h}-I)^{\ell}-\sum_{j=1}^{\ell} {\ell\choose j} 2^{-j}(\tau^{\widetilde{X}}_{-h}-I)^{\ell+j}.$$

    Now start with $\ell=\Ga$. Since we are handling the inhomogeneous case, we have 
    $$\|(\tau^{\widetilde{X}}_{-h}-I)^{\Ga}g\|_{L^\infty}\leq 2^{-\Ga} \|(\tau^{\widetilde{X}}_{-2h}-I)^{\Ga}g\|_{L^\infty}+C\|g\|_{\widetilde{\GL'}^{\times}_{\Ga}}|h|^{\Ga}.$$
    Iterating until $2^N|h|\in [\frac{1}{2},1)$, we have 
    \begin{align*}
        \|(\tau^{\widetilde{X}}_{-h}-I)^{\Ga}g\|_{L^\infty}&\leq 2^{-N\Ga}\|(\tau^{\widetilde{X}}_{-2^Nh}-I)^{\Ga}g\|_{L^\infty}+C_{\Ga} N \|g\|_{\widetilde{\GL'}^{\times}_{\Ga}}|h|^{\Ga}\\&\leq 2^{-N\Ga}2^{\Ga}\|g\|_{\widetilde{\GL'}^{\times}_{\Ga}}+C_{\Ga}N\|g\|_{\widetilde{\GL'}^{\times}_{\Ga}}|h|^{\Ga} \leq C_{\Ga} \|g\|_{\widetilde{\GL'}^{\times}_{\Ga}}|h|^\Ga\left(1+\log\frac1{|h|}\right)
    \end{align*}

    Next, we consider $\ell=\Ga-1$, then if $2^{J}|h|\in [\frac{1}{2},1)$,
    \begin{align*}
        \|(\tau^{\widetilde{X}}_{-h}-I)^{\Ga-1}g\|_{L^\infty}&\leq 2^{-\Ga+1} \|(\tau^{\widetilde{X}}_{-2h}-I)^{\Ga-1}g\|_{L^\infty}+C'\|g\|_{\widetilde{\GL'}^{\times}_{\Ga}}|h|^{\Ga}\left(1+\log\frac1{|h|}\right)+C\|g\|_{\widetilde{\GL'}^{\times}_{\Ga}}|h|^{\Ga}\\
        &\leq \cdots \\
        &\lesssim 2^{-J(\Ga-1)}\|g\|_{\widetilde{\GL'}^{\times}_{\Ga}}+\|g\|_{\widetilde{\GL'}^{\times}_{\Ga}}|h|^{\Ga}\sum_{k=0}^{J-1}2^{k(\Ga-\ell)}\bigg(1+\log\frac{1}{2^k|h|}\bigg) \\
        %&\leq C'|h|^{\Ga-1}\|g\|_{\widetilde{\GL'}^{\times}_{\Ga}}+C'\|g\|_{\widetilde{\GL'}^{\times}_{\Ga}}|h|^{\Ga}\bigg[\left(1+\log\frac1{|h|}\right)^2+\left(1+\log\frac1{|h|}\right)\bigg]\\
        &\leq C''|h|^{\Ga-1}\|g\|_{\widetilde{\GL'}^{\times}_{\Ga}}.
    \end{align*}
    Here, note that 
    $$\sum_{k=0}^{J-1}2^{k(\Ga-\ell)}\bigg(1+\log\frac{1}{2^k|h|}\bigg)=2^J\sum_{j=1}^{J}2^{-j}\log\bigg(\frac{e2^j}{2^J|h|}\bigg)\lesssim2^J\sum_{j=1}^J2^{-j}(1+j)\lesssim 2^J\simeq |h|^{-1}$$
    Inductively, we have the desired estimate.
\end{proof}

\begin{lemma}\label{lemma2}
    For fixed $R>1$, we have 
    $$\sup_{c\in \C^n}\|F_c\|_{\GL_{\Ga}(B(0,R))}\leq C_{\Ga,R}\|g\|_{\widetilde{\GL'}^{\times}_{\Ga}}.$$
\end{lemma}
\begin{proof}
    Let $0<|h|<1$ and $u,u+(\Ga+1)h\in B(0,R)$. Using Lemma \ref{lemma1} and Marchaud's inequality, we have 
    \begin{align*}
        |\Delta^{\Ga+1}_hF_c(u)|&\leq \sum_{\ell=0}^{\Ga-1}{\Ga+1\choose \ell}|h|^{\Ga+1-\ell}C'\|g\|_{\widetilde{\GL'}^{\times}_{\Ga}}|h|^{\ell}+\Ga C_{\Ga}|h|\|g\|_{\widetilde{\GL'}^{\times}_{\Ga}}|h|^{\Ga}\left(1+\log\frac{1}{|h|}\right)+\|g\|_{\widetilde{\GL'}^{\times}_{\Ga}} |h|^{\Ga} \\
        &\leq C_{\Ga, R} \|g\|_{\widetilde{\GL'}^{\times}_{\Ga}} |h|^{\Ga}.
    \end{align*}
    We then obtain the desired estimate from the fact that 
    $$\|F_c\|_{L^\infty(B(0,R))}=\|g\|_{L^\infty(B(c,R))}\leq \|g\|_{L^\infty(\C^n)}\leq\|g\|_{\widetilde{\GL'}^{\times}_{\Ga}}.$$
   
    Here we understand $\GL_{\Ga}$ as the finite-difference Zygmund space, and it is known that $\GL_{\Ga}(B(0,2))\hookrightarrow C^{\Ga-1}(B(0,2))$.
\end{proof}

Next, we use the equivalence between (local) classical Campanato space and Zygmund space to obtain the desired result. More precisely, we want to show that
\begin{align*}
\sup_{0<r<1}
r^{-\Ga}
\inf_{P\in\mathscr P_\Ga}
\|f-P\|_{L^\infty(B(0,r))}
\lesssim
\|f\|_{\Lambda_{\Ga}(B(0,2))}.
\end{align*}

In fact, by the classical Whitney inequality for polynomial approximation on balls (see \cites{Bru70, BK00}), the local modulus of smoothness of order $\Ga+1$ controls the best uniform approximation by polynomials of total degree at most $\Ga$. More precisely, for fixed $r\in (0,1)$,
\begin{align*}
    \inf_{P\in\mathscr P_\Ga}
\|f-P\|_{L^\infty(B(0,r))}&\lesssim \sup_{z,z+(\alpha+1)h\in \ol{B(0,r)}}|\Delta_h^{\Ga+1}f(z)| \\
&\lesssim_{\Ga} r^{\Ga}\sup_{z,z+(\Ga+1)h\in \ol{B(0,r)}}\frac{|\Delta_h^{\Ga+1}f(z)|}{|h|^{\Ga}}\\
&\lesssim r^{\Ga} \sup_{|h|<\frac{2r}{\Ga+1}}\frac{\|\Delta_h^{\Ga+1}f\|_{L^\infty(B(0,2))}}{|h|^{\Ga}}\lesssim r^{\Ga}\|f\|_{\GL_{\Ga}(B(0,2))}
\end{align*}
and the implicit constant is independent of $r$; therefore,
\begin{align*}
\sup_{0<r<1}
r^{-\Ga}
\inf_{P\in\mathscr P_\Ga}
\|f-P\|_{L^\infty(B(0,r))}
\lesssim
\|f\|_{\GL_{\Ga}(B(0,2))}.
\end{align*}

Meanwhile, 
\begin{align*}
&    \sup_{0<r<1, c\in \C^n} r^{-\Ga}\inf_{P\in\mathscr P_\Ga}
\|\tau^{X}_cg-P\|_{L^2(B(0,r);dx/r^{2n})}\\&=\sup_{0<r<1, c\in \C^n}
r^{-\Ga}
\inf_{P\in\mathscr P_\Ga}
\|g(\cdot+c)\Go(z,c)-P\|_{L^2(B(0,r);dx/r^{2n})} \\
&=\sup_{0<r<1, c\in \C^n}
r^{-\Ga}
\inf_{P\in\mathscr P_\Ga}
\|g(\cdot)-P(\cdot-c)\Go(c,\cdot)\|_{L^2(B(c,r);dx/r^{2n})}.
\end{align*}
Combining with Lemma \ref{lemma2}, we have
\begin{align*}
&\sup_{0<r<1, c\in \C^n}r^{-\Ga}
\inf_{P\in\mathscr P_\Ga}
\|g(\cdot)-P(\cdot-c)\Go(c,\cdot)\|_{L^2(B(c,r);dx/r^{2n})}\\
&=  \sup_{0<r<1, c\in \C^n} r^{-\Ga}\inf_{P\in\mathscr P_\Ga}
\|F_c-P\|_{L^2(B(0,r);dx/r^{2n})} \\
&\lesssim \sup_{ c\in \C^n}\|F_c\|_{\GL^{\Ga}(B(0,2))}\lesssim C_{\Ga}\|g\|_{\widetilde{\GL'}^{\times}_{\Ga}}.
\end{align*}

Finally, since $g\in L^\infty(\C^n)$, we have 
$$\sup_{B:\ r(B)\geq 1} \bigg(\frac{1}{r^{2n}}\int_B|g(z)|^2dz\bigg)^{\frac{1}{2}}\lesssim \|g\|_{ L^\infty}\lesssim \|g\|_{\widetilde{\GL'}^{\times}_{\Ga}}.$$
Therefore, $\widetilde{\GL'}^{\times}_{\Ga}\hookrightarrow\GL^{\times}_{\Ga}$.

Combining with Theorem \ref{e1} and Lemma \ref{e2}, we conclude the main theorem of this section.
\begin{theorem}
    For $\Ga>0$, the spaces $\GL^\times_{\Ga}\cong \widetilde{\GL}^{\times}_{\Ga}$ with equivalent norm. 
\end{theorem}

\begin{remark}
    The equivalence between $\widetilde{\GL}^{\times}_{\Ga}$ and $\widetilde{\GL'}^{\times}_{\Ga}$ is also true for $\Ga\notin \N$ with  the definition of 
    $$\|f\|_{\widetilde{\GL'}^{\times}_{\Ga}}:=\|f\|_{L^\infty}+\sup_{|h|<1,z\in\C^n}\frac{|(\tau^{\widetilde{X}_h}-I)^{\lfloor\Ga\rfloor+1}f(z)|}{|h|^{\Ga}}.$$ 
    In this case, the logarithmic term in the Marchaud's inequality does not occur. Thus, for $\Ga>0$, 
    $$\GL^{\times}_{\Ga}\cong\widetilde{\GL}^{\times}_{\Ga}\cong\widetilde{\GL'}^{\times}_{\Ga}$$
    with equivalent norms.
\end{remark}

\section*{Acknowledgement}
We would like to thank Riju Basak for introducing this topic to us.

Use of AI: During the preparation of this manuscript, the author used ChatGPT (OpenAI) to improve the language, assist in the search for relevant literature, and discuss and check certain mathematical arguments. The author has reviewed the resulting content and takes full responsibility for the mathematical results, proofs, and final manuscript.

\end{document}